\documentstyle{elsart}

\font\nineit=cmti12
\font\ninebf=cmbx12

\def\ifm#1#2{\relax\ifmmode#1\else#2\fi} 

\newcommand{\xon}    {\ifm {X_1,\ldots,X_n} {$X_1,\ldots,X_n$}}
\newcommand{\xo}[1]  {\ifm {X_1,\ldots,X_{#1}} {$X_1,\ldots,X_{#1}$}}
\newcommand{\yon}    {\ifm {Y_1,\ldots,Y_n} {$Y_1,\ldots,Y_n$}}

\newcommand{\fon}    {\ifm {f_1,\ldots,f_n} {$f_1,\ldots,f_n$}}
\newcommand{\ifon}   {\ifm {(f_1,\ldots,f_n)} {$(f_1,\ldots,f_n)$}}
\newcommand{\ifo}[1] {\ifm {(f_1,\ldots,f_{#1})} {$(f_1,\ldots,f_{#1})$}}

\newcommand{\fo}[1]  {\ifm {f_1,\ldots,f_{#1}} {$f_1,\ldots,f_{#1}$}}
\newcommand{\Zxon}   {\ifm {\Z[X_1,\ldots,X_n]} {$\Z[X_1,\ldots,X_n]$}}
\newcommand{\Rxon}   {\ifm {R[X_1,\ldots,X_n]} {$R[X_1,\ldots,X_n]$}}
\newcommand{\Ryon}   {\ifm {R[Y_1,\ldots,Y_n]} {$R[Y_1,\ldots,Y_n]$}}
\newcommand{\Kxon}   {\ifm {K[X_1,\ldots,X_n]} {$K[X_1,\ldots,X_n]$}}
\newcommand{\Bxon}   {\ifm {B[X_1,\ldots,X_n]} {$B[X_1,\ldots,X_n]$}}
\newcommand{\Byon}   {\ifm {B[Y_1,\ldots,Y_n]} {$B[Y_1,\ldots,Y_n]$}}
\newcommand{\Qxon}   {\ifm {\Q[X_1,\ldots,X_n]} {$\Q[X_1,\ldots,X_n]$}}
\newcommand{\Qxo}[1]   {\ifm {\Q[X_1,\ldots,X_{#1}]} {$\Q[X_1,\ldots,X_{#1}]$}}
\newcommand{\Qyon}   {\ifm {\Q[Y_1,\ldots,Y_n]} {$\Q[Y_1,\ldots,Y_n]$}}
\newcommand{\klk}    {\ifm {,\ldots,} {$,\ldots,$}}
\newcommand{\plp}    {\ifm {+\cdots+} {$+\ldots+$}}

\newcommand{\bino}[2] {\ifm {{#1}\choose{#2}} {${#1}\choose{#2}$}}
\newcommand{\abs}[1]  {\ifm {\vert{#1}\vert} {$\vert{#1}\vert$}} 
\newcommand{\norm}[1] {\ifm {\vert\vert{#1}\vert\vert} {$\vert\vert{#1}\vert\vert$}} 

\newcommand{\Tr}  {\mbox{\it Tr}}
\newcommand{\Trt} {\ifm {\widetilde{\mbox{\it Tr}}} {$\widetilde{\mbox{\it Tr}}$} }

\newcommand{\spar} {\vskip 0.0cm}

\def \R{{\rm I\kern -2.2pt R\hskip 1pt}} 
\def \Z{{\rm Z\!\!Z}}
\def \N{{\rm I\kern -2.1pt N\hskip 1pt}} 

\def \K {{\it K}}
\def \P{{\rm I\kern -2.2pt P\hskip 1pt}} 

\newcommand{\Strichq}  {\vrule height0.65em width0.05em depth0em \,}
\newcommand{\Q}        {\ifm {\mbox{\rm Q}\hspace{-0.54em}\Strichq\>\;}
                             {$\mbox{\rm Q}\hspace{-0.54em}\Strichq\>\;$}}
\newcommand{\C}        {\ifm {\mbox{\rm C}\hspace{-0.45em}\Strichq\>}
                             {$\mbox{\rm C}\hspace{-0.45em}\Strichq\>$}}

\newtheorem{notation}{Notation} 
\newtheorem{remark}{Remark}

\begin{document}


\begin{frontmatter}

\title{LOWER BOUNDS FOR DIOPHANTINE APPROXIMATIONS\thanksref{funds}} 

\author{M. Giusti}
\address{GAGE, Centre de Math\'ematiques. \'Ecole Polytechnique. F-91228
Palaiseau Cedex. France}

\author{J. Heintz, K. H\"agele, J. E. Morais, L. M. Pardo} 
\address{Departamento de Matem\'aticas, Estad\'{\i}stica y Computaci\'on,
Facultad de Ciencias, Universidad de Cantabria, E-39071 Santander, Spain}

\author{J. L. Monta\~na}
\address{Departamento de Matem\'atica e Inform\'atica, Campus de
Arrosad\'{\i}a, Universidad P\'ublica de Navarra, E-31006 Pamplona, Spain}

\thanks[funds]{Research was partially supported by the following French and
Spanish grants~: GDR CNRS 1026 MEDICIS, DGCYT PB92--0498--C02--01 and
PB93--0472--C02--02}


\begin{abstract} 
We introduce a subexponential algorithm for geometric solving of
multivariate polynomial equation systems whose bit complexity depends mainly
on intrinsic geometric invariants of the solution set.
>From this algorithm, we derive a new procedure for the decision of
consistency of polynomial equation systems whose bit complexity is
subexponential, too. As a byproduct, we analyze the division of a polynomial
modulo a reduced complete intersection ideal and from this, we obtain an
intrinsic lower bound for the logarithmic height of diophantine
approximations to a given solution of a zero--dimensional polynomial
equation system. This result represents a multivariate version of
Liouville's classical theorem on approximation of algebraic numbers by
rationals.
A special feature of our procedures is their {\em polynomial\/} character
with respect to the mentioned geometric invariants when instead of bit
operations only arithmetic operations are counted at unit cost.
Technically our paper relies on the use of straight--line programs as a data
structure for the encoding of polynomials, on a new symbolic application of
Newton's algorithm to the Implicit Function Theorem and on a special, basis
independent trace formula for affine Gorenstein algebras.
\end{abstract}

\end{frontmatter}


\typeout{Section 1}
\begin{section}{Introduction}\label{intro}

The present pages represent a continuation of
\cite{gihemomorpar,gihemopar,pardo}. These papers concern the design of
algorithms of intrinsic type to solve systems of polynomial equations. 
Solving is then applied to decide consistency of systems of polynomial
equations. By ``intrinsic type" we mean algorithms which are able to
distinguish between the semantical and the syntactical character of the
input system in order to profit from both for the improvement of the
complexity estimates.  \spar

With respect to bit complexity we show how the time necessary to solve a
given polynomial equation system is related to the affine degree and the
(affine) logarithmic height of the corresponding diophantine variety. In
this sense, the results of \cite{gihemomorpar} and \cite{gihemopar} show
already that the (affine) degree of an input system is associated with the
complexity when measured in terms of number of arithmetic operations.
However, there are still some drawbacks in this approach.  \spar

The first one is that the algorithms developed in \cite{gihemopar} require
iterative calls to the procedure and that the algorithms in
\cite{gihemomorpar} rely on the use of algebraic numbers. Thus in the case
of \cite{gihemomorpar} the algorithms are not ``rational" although their
inputs and outputs are. \spar

A second drawback concerns the modeling used to measure complexity.  The
quantity of arithmetic operations of an algorithm does not explain
sufficiently what happens when we are using it on a ``real life" computer.
Years of experience show that models of bit complexity (Turing machines,
random access machines or equivalent models) represent more realistic
patterns for practical computing. In this sense a study of bit complexity of
the intrinsic complexity of the algorithms of \cite{gihemomorpar} and
\cite{gihemopar} becomes necessary. \spar

In the present paper we deal with both disadvantages of the algorithms in
\cite{gihemomorpar} and \cite{gihemopar}, giving practicable solutions (cf.
also Section \ref{algoritmo} below).  First, our new algorithm is completely
rational.  It does not require any constants other than those in the field
of coefficients of the input system. Secondly, to improve the bit complexity
estimates, we introduce a suitable notion of height of affine diophantine
varieties which is inspired by the corresponding notion introduced for
projective varieties in \cite{BeYg3,bgs2,Faltings,Nesterenko1,Nesterenko2,philippon0,philippon1,philippon2,philippon3}. As shown in
Section 2, our notion of height is strongly related to the bit complexity of
geometric elimination procedures of any kind. \spar

Our notion of height combined with a new algorithmic interpretation of
duality theory for complete intersection ideals yields a new Liouville
estimate (cf. Sections \ref{Liouville} and \ref{division} below).  Liouville
estimates can be applied to get lower time bounds for the numerical analysis
approach to solving systems of polynomial equations. In particular, we show
the practical inefficiency of both floating point and binary encoding for
rational approximation of zero--dimensional multivariate polynomial equation
systems even in the algorithmically well suited cases. \spar

To illustrate these introductory observations, let us first consider the
following two problems:

\begin{prob}\label{problem1}

Let be given a sequence of $n$ integer polynomials of degree at
most $3$, with coefficients in $\{0,1\}$, concerning three variables each:

$$p_1(X_{1,1},X_{1,2},X_{1,3}) \klk p_n(X_{n,1},X_{n,2},X_{n,3})$$ 

Decide whether the following system of polynomial equations
has a solution:

$$X_1^2-X_1=0 \klk X_n^2-X_n=0, Y_1^2-Y_1=0 \klk Y_{3n}^2-Y_{3n}=0$$

$$\prod_{i=1}^np_i(X_{i,1},X_{i,2},X_{i,3})-(1+ \sum_{j=1}^{3n} Y_j2^j)=0$$
\end{prob}

\begin{prob}\label{problem2}

Let be given in binary representation an integer $k\in\N$ and the polynomials: 

$$X_1^2-X_1=0 \klk X_n^2-X_n=0 $$ 

$$k - (X_1+2X_2 \plp 2^{n-1}X_n)=0$$ 

Decide whether this system has a solution. 
\end{prob}

Both problems have a similar ``syntactical form" (i.e. they look very
similar, as a consistency question for ``syntactically easy" polynomial
equation systems). However, they possess completely different ``semantical"
characters~: the first problem is a translation of a well--known {\bf
NP}--complete problem (3-satisfiability, for short 3SAT), while the second
one just concerns the binary encoding of $k$ with $n$ bits (if there exists
such an encoding). This means that the second system is consistent if and
only if $k$ may be encoded using at most $n$ binary digits (cf. \cite{HeMo,pardo}). Traditional symbolic procedures deal with both problems using
in each case the same generalistic treatment.  However these two problems
demand for different algorithms that may profit from their different
semantical features.  The construction of algorithms which are able to
distinguish between equation systems which are semantically well suited and
such which are not is the main goal of our paper. We shall also consider the
question of consistency of polynomial equation systems in the following
terms:

\begin{prob}[Effective Nullstellensatz]\label{problem3} 
Given a sequence of polynomials $f_1\klk f_s \in \Z[X_1 \klk
X_n]$,  decide whether the affine algebraic variety 

$$V(f_1 \klk f_s):=\{x\in \C^n: f_1(x)=\cdots=f_s(x)=0\}$$ 

is empty or not.
\end{prob}

As both Problems \ref{problem1} and \ref{problem2} above can be written as a
special case of this more general Problem \ref{problem3}, we have to find a
way to distinguish (within the context of Problem \ref{problem3}) between
the different levels of difficulty Problems \ref{problem1} and \ref{problem2}
represent. To design algorithms which solve this problem, taking
care of the special features of the particular instance of the input system,
we need two major geometric invariants: the affine degree and the affine
height of the system (cf. Section \ref{intrinsic} below).  Assuming these
two notions we are going to show in this paper the following result:

\begin{thm}\label{teouno} There exists a bounded error probability
Turing machine which solves the following task: given polynomials $f_1 \klk
f_{n+1}\in\Zxon$ of degree at most 2 and of (logarithmic) height $h$, the
Turing machine decides whether the variety $V\ifo {n+1}$ is empty or not.
Moreover, if $\delta$ is the intrinsic (affine) degree of the system and
$\eta$ is the intrinsic (logarithmic) height of the system, the Turing
machine answers in (bit) time

$$(n h \delta\eta)^{O(log_2n)}$$ 

using a total amount of arithmetic operations in \Q of

$$(n\delta)^{O(1)}.$$ 

\end{thm}

Our algorithm first solves a suitable polynomial equation system and then
uses this information for the consistency test of the original system $\fo
{n+1}$. Solving is done inductively. In order to explain the procedure let
us assume that the polynomials $\fon$ form a regular sequence, each ideal
$\ifo i$ $1\leq i \leq n$, being radical. Then the algorithm proceeds in $n$
steps, solving at each stage $1\leq i \leq n$ the system $f_1=
0,\ldots,f_i=0$.
\spar

The corresponding intermediate algebraic varieties are obtained by a lifting
process from special zero-dimensional varieties which we call the {\sl
lifting fibers}.  The lifting process is based on a division-free symbolic
version of the Newton-Hensel algorithm and represents a algorithmic version
of the Implicit Function Theorem. In this way, the procedure constructs the
lifting fiber of step $i+1$ from the lifting fiber of step $i$. Each
inductive step includes two cleaning phases~: one is performed to throw away
extraneous projective components (mainly those at infinity) and a second one
reduces the size of the representation of the integers which appear during
the process. This second cleaning phase is included to avoid uncontrolled
growth of integer coefficients of intermediate polynomials. \spar

Finally we show how the following two computational problems are related:
the division problem in the Nullstellensatz and the numerical analysis
approach to solving systems of multivariate polynomial equations.

\begin{prob}[Division Problem in the Effective Nullstellensatz]
\label{problem4} Let\\ be given polynomials
$f_1,\ldots,f_{s}\in\Z[X_1,\ldots,X_n]$ without any common zeroes (i.e. the
polynomials satisfy the condition $V(f_1 \klk f_{s}) = \emptyset$), compute
a non-zero integer $a\in\Z$, and polynomials
$g_1,\ldots,g_s\in\Z[X_1,\ldots,X_n]$ such that

$$a = g_1 f_1 \plp g_s f_s$$ 

holds.
\end{prob}

\begin{prob}[Numerical Analysis Approach to Solving]\label{problem5}
Given a real \linebreak 
number $\varepsilon >0$ and a regular sequence of polynomials
$\fon$$\in\Zxon$ of degree at most $d$ defining a zero-dimensional affine
variety \linebreak
$V:= V\ifon$, compute for any point $\alpha\in V$ an approximation
$a\in\C^n$ at level $\varepsilon>0$, i.e. find by an effective procedure a
point $a\in\C$ such that

$$\norm{a-\alpha} < \varepsilon$$ 

holds.
\end{prob}

At first glance Problems \ref{problem4} and \ref{problem5} seem to be
unrelated. However we shall see how the solution of Problem \ref{problem4}
provides lower bounds for any solution of Problem \ref{problem5}. We shall
state this observation in terms of Liouville estimates. \spar

In order to understand the relation between Problems \ref{problem4} and
\ref{problem5}, let us observe that in the context of Problem \ref{problem5}
any approximation a computer may output is necessarily rational (i.e. such
an approximation must be a point $a$ belonging to $\Q[i]^n$). Such an output
is assumed to be encoded in binary (i.e. by the binary expansion of
denominators and numerators of the coordinates of $a$). \spar

In the past the division problem in the Nullstellensatz was studied in both
number theory and computer science simultaneously. The first results mainly
established upper bounds for the degree
(\cite{brown1,can-gal-he,cagahe,kollar}) and height (\cite{BeYg1,BeYg2}) of
the polynomials $g_1\klk g_s$ appearing in Problem \ref{problem4}. The
complexity study in (\cite{figismi,krick-pardo1,krick-pardo:CRAS}) yields
optimal upper bounds for both invariants (see also \cite{Elkadi}) and can be
applied to derive Liouville estimates in the following sense (cf. also
Subsection \ref{step} below): \spar

A Liouville estimate is a lower bound for the binary length (the height) of
the numerators and denominators of the coordinates of $a$ in terms of the
variety $V$, the point $\alpha$ and the real number $\varepsilon$ in the
statement of Problem \ref{problem5} above. The following Theorem
\ref{Liouville} gives such a Liouville estimate. \spar

Let us introduce the following notions and notations: let $V\subset \C^n$ be
a zero-dimensional diophantine variety not intersecting $\Q[i]^n$ and let
$\alpha$ be a point of $V$. Let $\varepsilon>0$ be a real number. We call a
point $a=({{a_1}\over{q_1}}, \ldots, {{a_n}\over{q_n}}) \in \Q[i]^n$ a
rational approximation of $\alpha$ at level $\varepsilon>0$, if

$$\norm{\alpha - a} < \varepsilon$$

holds. (Here \norm{\cdot} denotes the norm associated with the usual
hermitian product in $\C^n$).

A regular sequence $f_1,\ldots,f_n\in\Z[X_1,\ldots,X_n]$ is said to be
smooth if for any $1\leq i \leq n$ the Jacobian $J(f_1,\ldots,f_i)$ is not a
zero-divisor modulo the ideal $(f_1,\ldots,f_i)$. Observe that for a smooth
regular sequence the ideals $\ifo i$ are always radical.

\begin{thm}\label{cotainf}\label{Liouville} There exists a universal
constant $C$ with the following property: given a smooth regular
sequence of polynomials $f_1,\ldots,f_n\in \Z[X_1,\ldots,X_n]$ of degree at
most $d$ and of logarithmic height at most $h$. Furthermore let be given a
point $a\in\Q[i]^n$ and a real number $0 < \varepsilon \leq 1$. Suppose that
the variety $V:=V(f_1,\ldots,f_n)$ verifies the following conditions~:

\begin{itemize}
\item $V\cap \Q[i]^{n}=\emptyset$.
\item there exists a point $\alpha \in V$ such that 

$$\norm {a-\alpha} < \varepsilon \leq 1$$
\end{itemize}

holds. Then we have for any denominator $q$ of $a$ the following
inequality~:
 
$$ {{-log_2 \varepsilon}\over {(nd\delta)^{C}}}-(h+\eta)\leq
log_2 q,$$

where $\delta$ denotes the degree and $\eta$ the (affine) logarithmic height
of $V$.  \end{thm}

\begin{pf} According to Proposition \ref{estrategia} below we have:

$${{\varepsilon^{-1}} \over { \abs {g_{n+1}(\alpha)} (2\norm \alpha +1) }}
\leq q^2.$$

Thus, we look for upper bounds for $\norm \alpha$ and $\abs
{g_{n+1}(\alpha)}$.
First, from Proposition \ref{ht-cota} below we have:

$$(2 \norm \alpha +1 ) \leq \sqrt n \delta 2^{2\eta}.$$

On the other hand, combining Theorem \ref{paso-division} and Lemma
\ref{cotas} below we conclude:

$$log_2\abs{g_{n+1}(\alpha)} \leq (n d \delta)^{O(1)}(h+\eta+log_2 q).$$

>From these two upper bounds we easily deduce the bounds of the theorem.
\qed
\end{pf}

We observe that for $n=1$ Theorem \ref{Liouville} represents just a
restatement of Liouville's classical result with somewhat coarser bounds
(compare e.g. \cite[Kapitel I, Satz 1]{Schneider}).
\spar

The height estimates in \cite{BeYg1,BeYg2,Elkadi,krick-pardo1} and
\cite{krick-pardo:CRAS}, combined with the methods described in subsection
\ref{step} below (Proposition \ref{estrategia}) produce Liouville bounds
that relate the syntactical description of $V$ as given by the input, the
approximation $\varepsilon$ and the height of the denominator of $a$. In
fact, these estimates yield the inequality~:

\begin{equation}\label{ugl1}
 -log_2\varepsilon \leq d^{O(n)} h log_2q
\end{equation}

On the other hand: degree and height of the solution set represent
reasonable lower time bounds for the numerical analysis approach to
polynomial equation solving by rational point approximation. Exponential
degrees and exponential height produce exponential lower time bounds for
numerical methods of polynomial equation solving (cf. \cite{pardo}). \spar

A lower bound such as inequality (\ref{ugl1})  above, i.e. a lower bound for
$log_2 q$, represents a lower bound for the output length and therefore a
lower time bound for numerical methods of polynomial system solving.  The
results of \cite{BeCaRoSz,BeYg1,BeYg2,bgs1,bgs2,Elkadi,krick-pardo1,krick-pardo:CRAS,philippon1,philippon2,philippon3} 
imply that an
approximation level of $\varepsilon:=2^{d^{O(n)}}$ is sufficient in order to
characterize (and to distinguish adequately) the solutions of the variety
$V$ in Problem \ref{problem4}. On the other hand the examples of
\cite[Chapitre 4, Proposition 12]{Mignotte} and \cite{Slietal,Slietal2} show
that even in case of semantically and syntactically ``easy" systems such an
approximation level may be necessary. Therefore it becomes reasonable to fix
an approximation level for numerical solving of $\varepsilon:=2^{d^n}$,
where $d^n$ represents the B\'ezout number of the input system. Under this
assumption, inequality $(1)$ becomes meaningless and more precise lower
bounds are required. \spar
 
However an exponential lower time bound for numerical solving is implicitly
contained in the following corollary to Theorem \ref{Liouville}:

\begin{cor}\label{Cor-Liouville} Given a smooth regular sequence
$f_1,\ldots,f_n\in\Z[X_1,\ldots,X_n]$, and points $\alpha\in V$,
$a\in\Q[i]^n$ verifying the conditions of Theorem \ref{Liouville} above, let
$\varepsilon >0$ be a level of approximation such that $-log_2 \varepsilon =
d^n$ (where $d^n$ represents the B\'ezout number of the system). Then we
have

$$ {{d^n}\over {(nd\delta)^{C}}}-(h+\eta)\leq log_2 q,$$ 

where $C$ is a suitable universal constant as in Theorem \ref{Liouville}. In
particular, for systems of ``small" degree and height, the output length for
numerical solving methods is necessarily exponential in the number of
variables. 
\end{cor}

For instance, we consider the following smooth regular sequence of
quadratic polynomials:

$$X_1^2+X_1+1,X_2-X_1^2,\ldots,X_n-X_{n-1}^2$$

This sequence defines a zero-dimensional variety with just two points (the
degree of the variety is 2) of small height (the height is 1). Thus, for a
level of approximation $\varepsilon>0$ with $-log_2\varepsilon = 2^n$, the
lower bound obtained from inequality (\ref{ugl1}) says just:

$$2^n\leq 2^{O(n)} log_2 q,$$

whereas the lower bound from Corollary \ref{Cor-Liouville} states that every
rational approximation of level $\varepsilon$ of either solution of this
system has exponential binary length. Thus, the main consequence of this
corollary is that both binary and floating point encoding of numbers in
$\Q[i]$ are not suitable to reach the appropriate level of approximation. An
alternative encoding is therefore required. \spar

Another way out of this dilemma may consist in the following approach
initiated in \cite{ShSm93a} and further developed in \cite{ShSm93b,ShSm93c,ShSm93d,ShSm1,Dedieu1,Dedieu2} instead of approximating rationally (or
similarly by floating point arithmetic) the solutions of the
zero--dimensional input system \fon up to the appropriate approximation
level $\varepsilon = 2^{d^n}$, we just try to find approximate zeroes (in
the sense of \cite{ShSm1}) of the system, i.e. we try to find points
$a\in\Q[i]^n$ from which a suitable version of Newton's algorithm converges
quadratically to a true zero of the system.

\end{section} 

\typeout{Section 2}

\begin{section}{Notions, notations and Results}\label{notaciones}

We are going to study the consistency question of Problem \ref{problem4}
(i.e. the decisional problem in the effective Nullstellensatz) only when the
input system of polynomial equations $f_1,\ldots,f_r\in\Z[X_1,\ldots,X_n]$
verifies the conditions:

\begin{itemize}
\item $r\leq n+1$,
\item the sequence $f_1\klk f_{r-1}$ is a smooth regular sequence 
\end{itemize}

These two additional conditions are not really restrictive. As shown in
\linebreak
\cite{figismi,giu-he92,gihesa,krick-pardo1,krick-pardo:CRAS} an
effective version of Bertini's Theorem allows to reduce a general input
system $f_1,\ldots,f_s\in \Z[X_1,\ldots,X_n]$ to such a system. We just have
to find generic $\Z$-linear combinations of the input polynomials such that
these linear combinations contain only coefficients in $\Z$ of logarithmic
height $O(n\cdot log_2 d)$, where $n$ is the number of variables and $d$ is
an upper bound for the degrees of the input polynomials. Under the
assumption of this preprocessing we may suppose without loss of generality
that our input equations satisfy the following conditions:

\begin{itemize} 
\item the ideals $(f_1,\ldots,f_i)$ are radical, $1\leq i\leq r-1$, 
\item for $1\leq i\leq r-1$ the varieties $V_i = V(f_1,\ldots,f_i)$ are complete
intersection affine varieties of dimension $n-i$.
\end{itemize}

\begin{notation} Given $R\subset B$ an extension of a commutative ring
(which converts $B$ into a $R$-algebra), and an element $b\in B$, we denote
by $\eta_b:B \longrightarrow B$ the $R$-linear endomorphism induced by the
multiplication of the elements of $B$ by $b$ (in the following we shall call
such a linear map a homothety). If $B$ is a free $R$-module of finite rank,
we denote by $M_b$ the matrix of the homothety $\eta_b$ and by $\chi_b\in
R[T]$ the characteristic polynomial of $\eta_b$. Moreover, if $R$ is a
unique factorization domain, we shall denote by $m_b$ the primitive minimal
polynomial of $\eta_b$. Observe that $\chi_b$ and $m_b$ are monic
polynomials of $R[T]$, which for short we will call characteristic and
minimal polynomial of $b$, respectively. \end{notation}

In order to decide whether $(f_1,\ldots,f_r)$ represents the trivial ideal
we just need to compute the following items:

\begin{itemize} 

\item A linear change of coordinates $(X_1,\ldots,X_n)\longrightarrow
(Y_1,\ldots,Y_n)$ such that the following $\Q$-algebra homomorphism represents
a Noether normalization of the variety $V(\fo {r-1})$, with the variables
$Y_1 \klk Y_{n-r+1}$ being free:

$$R:=\Q[Y_1,\ldots,Y_{n-r+1}] \longrightarrow
\Q[Y_1,\ldots,Y_n] / (f_1,\ldots,f_{r-1}) =: B.$$ 

Observe that this means that $B$ is an integral extension of $R$.

\item the matrix $M_{f_r}$ of the homothety $\eta_{f_r}$ with respect to an 
appropriate $R$-module basis.

\end{itemize} 

Suppose that such a Noether normalization is available. Then the $\Q$-algebra
homomorphism $R\rightarrow B$ is injective and $B$ is a free $R$-module of
rank at most $deg\ V_{r-1}$ (cf. \cite{gihesa,joos:tesis}).  Under
this assumption the ideal $(f_1,\ldots,f_r)$ is trivial if and only if the
matrix $M_{f_r}$ is unimodular what means that the determinant of $M_{f_r}$
is non zero and belongs to \Q. \spar

This comment shows how the original Problem \ref{problem3} of testing
consistency of polynomial equation systems can be reduced to the problem of
solving polynomial equation system in a very specific geometric sense. In
the next subsection we are going to explain what exactly we mean by this,
namely geometric solving. \spar


\begin{subsection}{Geometric Solving}\label{solving} 

The previous considerations reduce the search for a consistency test for the
polynomial equation system $f_1=0\klk f_n=0$, namely Problem \ref{problem3},
to the problem of computing a Noether normalization of the variety $V(\fo
{r-1})$ and the matrix $M_{f_r}$ of the homothety $\eta_{f_r}$ with respect
to a suitable $R$-module basis. Assume that \xon are already in Noether
position with respect to the variety $V(\fo {r-1})$, the variables $\xo
{n-r+1}$ being free. Then we consider the following integral ring
extension

$$R:=\Qxo {n-r+1} \longrightarrow \Qxon / \ifo {r-1} =:B .$$

Our assumptions on $\fo {r-1}$ imply that $B$ is a reduced algebra and a
finite free $R$-module. We are now going to explain what we mean by
``geometric solving" or ``geometric solution". First we need the following
notion of {\sl primitive element} of $B$:

\begin{defn}
Let $R$ be a ring of polynomials over $\Q$ and $R\subseteq B$ an integral
ring extension such that $B$ is reduced and $B$ is a free $R$-module of
finite rank. An element $u\in B$ is called a {\sl primitive element of the
ring extension $R\subseteq B$} if the degree of the minimal polynomial $m_u$
of $u$ equals the rank of $B$ as free $R$-module, i.e. if

$$deg\ m_u = rank_R\ B$$ 

holds.
\end{defn}

Let $\K$ be the quotient field of $R$ and $B'= \K\otimes_R B$ be the
localization of $B$ by the non-zero elements of $R$.  An element $u\in B$ is
a primitive element of $B$ if and only if for $D:=rank_r B = dim_k B'
$ the set $\{1,u \klk u^{D-1}\}$ represents a $\K$-vectorspace basis of
$B'$. \spar

The computation of the matrix $M_{f_r}$ is a consequence of the following
``generic point" description of the $\K$-algebra $B'$: 
this algebra is
characterized by the following items (which our algorithm will compute):

\begin{itemize} 

\item a $\K$-vectorspace basis of $B'$ 

\item for $n-r+2\leq i \leq n$ the matrices $M_{X_i}$ of the homotheties
$\eta_{X_i}: B \longrightarrow B$ with respect to the given basis (these
matrices describe the multiplication tensor of the \K-algebra $B'$ and hence
also of the $R$-algebra $B$)

\end{itemize}

We then obtain the matrix $M_{f_{r+1}}$ by substituting for the variables
\linebreak
$X_{n-r+2} \klk X_n$ appearing in the polynomial $f_{r+1}(\xon )$ the
matrices $M_{X_{x-r+2}} \klk M_{X_n}$. (In the sequel we shall write for
short \linebreak
$M_{f_{r+1}} = f_{r+1}(M_{X_{n-r+2}} \klk M_{X_n})$ for this
substitution, interpreting $f_{r+1}$ as an element of the polynomial ring
$R[X_{n-r+2} \klk X_n]$).
\spar

In this sense, geometric solving means just computing both a basis of the
\K-algebra $B'$ and the matrices $M_{X_{n-r+1}} \klk M_{X_n}$. This is done
making use of a suitable primitive element of $B$. \spar

In the context of this paper, the primitive element $u\in B$ will always be
chosen as the image in $B$ of a generic $\Z$-linear form of the variables
$X_{n-r+2}\klk X_n$. In particular, we may assume that $u$ is the image of a
linear form $U=\lambda_{n-r+2}X_{n-r+2} \plp \lambda_nX_n$, with
$\lambda_i\in\Z$ for $n-r+2\leq i\leq n$.  Let $T$ be a new indeterminate.
The minimal polynomial $m_u(T)$ of $u$ as an element of the $R$-algebra $B$
(or equivalently as an element of the $\K$-algebra $B'$) will be a monic
polynomial in $\Q[X_1\klk X_{n-r+1},T] = \Q \otimes_{\Z}R[T]$.  
\spar

We shall choose this minimal polynomial as an element of the ring
$\Z[X_1\klk X_{n-r+1},T]=R[T]$. In the sequel we shall pay special attention to the case
$r=n+1$, where we have $R=\Z$, $\K=\Q$ and where $m_u$ is a polynomial of
$\Z[T]$ of positive degree (the polynomial $m_u$ is then trivially monic
over $\Q[T]$ since we may divide it by its leading coefficient, which is a
non-zero integer). Discarding the content (the maximum common divisor) of
the coefficients of $m_u$, we may replace the minimal polynomial $m_u$ by its
primitive counterpart which we shall always denote by $q_u$. Similarly, for
the case $r\leq n$, we shall replace the minimal polynomial $m_u\in \Z[X_1
\klk X_{n-r+1},T]= R[T]$ by some ``cleaner" (however not necessarily
primitive) version $q_u$. \spar

Finally, as $\lbrace 1,u \klk u^{D-1}\rbrace $ is a basis of the
$\K$-vector space $B'$, for $n-r+2\leq i\leq n$ there are polynomials
$v_i^{(u)}\in R[T]$ and non-zero elements $\rho_i^{(u)}\in R$ such that
$\rho_i^{(u)}X_i-v_i^{(u)}(U)$ belongs to the ideal $\ifo {r-1}$ in
$\K[X_{n-r+2}\klk X_{n}]$. In particular, we have the following identity
between ideals of \linebreak
$\K[X_{n-r+2}\klk X_{n}]$:

$$\ifon = (q_u(U), \rho_1^{(u)}X_1-v_1^{(u)}(U) \klk
\rho_n^{(u)}X_n-v_n^{(u)}(U) ).$$

Moreover, if $M$ denotes the companion matrix of the homothety $\eta_u$ with
respect to the basis $\{1,u \klk u^{D-1}\}$, the matrices $M_{X_{n-r+2}}
\klk M_{X_n}$ characterizing the multiplication tensor of the $R$-algebra
$B$ (or equivalently of the $K$-algebra $B'$) are given by the formula

$$M_{X_i} = \rho_i^{(u)\ -1}\cdot v_i^{(u)}(M) $$

for $n-r+2\leq i\leq n$. To simplify notations we shall often omit the
superindex $(u)$ when refering to these polynomials. After these
explanations we shall assume without loss of generality that we have
$r:=n+1$ in the statement of Problems \ref{problem3} and \ref{problem4}.
Consequently we restrict the meaning of ``geometric solving" to the case
where we have given as input polynomials $\fon\in\Zxon$ forming a regular
sequence in \Qxon (we shall say in the future for short that $\fon\in\Zxon$
is a regular sequence). With these conventions ``geometric solving" means
the following:

\begin{defn}\label{geosolve}
An algorithm for geometric solving is a procedure which from a smooth
regular sequence $\fon\in\Zxon$ as input produces:

\begin{itemize}

\item a primitive element $u$ of the ring extension $\Q\longrightarrow
\Qxon/\ifon$ represented by a $\Z$-linear form $U=\lambda_1X_1 \plp
\lambda_nX_n$

\item the primitive minimal polynomial $q_u\in\Z[T]$ 

\item the parametrizations of the variety $V\ifon$ by the zeroes of $q_u$,
namely the (unique) primitive polynomials $ \rho_1^{(u)} X_1-v_1^{(u)}(T)
\klk \rho_n^{(u)} X_n-v_n^{(u)}(T) $ with $\rho_1^{(u)}\klk \rho_n^{(u)}$
non-zero integers and $v_1^{(u)}\klk v_n^{(u)}\in\Z[T]$ which satisfy the
conditions $max\{ deg\ v_1^{(u)}\klk deg\ v_n^{(u)}\} < deg\ q_u$ and $\ifon
= (q_u, \rho_1^{(u)} X_1-v_1^{(u)}(T) \klk \rho_n^{(u)} X_n-v_n^{(u)}(T) )$
in $\Qxon$. Moreover, the non-zero integer $\rho:=\prod_{i=1}^{n}\rho_i
\in\Z$, called a discriminant of the input system, is a multiple of the
discriminant of $q_u$.

\end{itemize}

\end{defn}

In the sequel we shall refer to the polynomials $U=\lambda_1X_1 \plp
\lambda_nX_n$ $\in\Zxon$ (with image $u$ in $\Qxon/\ifon$), $q_u(T)$,
\linebreak $v_1^{(u)}(T) \klk v_n^{(u)}(T) \in \Z[T]$ and $\rho_1^{(u)} \klk
\rho_n^{(u)} \in\Z$ as {\em geometric solution\/} of the equation system
$f_1=0\klk f_n=0$. \spar

Let us remark here that this notion of ``geometric solving'' has a long
history, which unfortunately we can not give here in full detail for mere
lack of space. One might consider \cite{ChGr} as an early reference where
this notion was implicitly used for the first time in complexity theory.

\end{subsection} 

\begin{subsection}{Intrinsic Parameters}\label{intrinsic} 

We said in the introduction that we are interested in algorithms for
geometric solving which are able to profit from ``good" geometrical
properties that an input system of polynomial equations might possess. This
makes it necessary to precise what such ``good" geometrical properties may
be and how to find measures for them. Therefore, we are going to define two
geometric invariants in this subsection that will arise as parameters of the
complexity of our procedures: the (affine) degree and the (affine)
logarithmic height of complete intersection varieties. The notion of degree
has been taken directly from \cite{joos:tesis}, while our notion of height
is strongly inspired by the corresponding notion developed for projective
varieties in \cite{Nesterenko1,Nesterenko2,philippon0,philippon1,philippon2,philippon3}. \spar

Let us first recall the notion of degree of an affine algebraic variety. We
begin by defining the degree of a zero-dimensional variety and then we
extend this notion to positive dimensional complete intersection algebraic
varieties. \spar

To fix notations, let $V\subseteq\C^n$ be an algebraic subset (variety)
given as the set of common zeroes in $\C^n$ of a smooth regular sequence of
polynomials $\fo i \in \Zxon$ of degree at most $d$. If $V$ is
zero-dimensional (i.e. if $i=n$) the degree of $V$ is defined as the number
of points of $V$ (points at infinity are not counted in this definition). 
In the general case (when $dim\ V = n-i\geq 1$ holds) let us consider the
class $\cal D$ of all affine linear subspaces of $\C^n$ of dimension $i$
(defined as the set of solutions in $\C^n$ of a linear equation system
$L_1=0 \klk L_{n-i} = 0$ where $L_k = a_{k1} X_1 \plp a_{kn} X_n + a_{k0}$
is an affine linear polynomial with coefficients $a_{kj}\in \Z$ for $1\leq
k\leq n-i$ and $ 0\leq j\leq n$). Let ${\cal D}_V$ be the subclass of $\cal
D$ of all these affine linear spaces $H\in \cal D$ such that $H\cap V$ is a
zero-dimensional variety.

\begin{defn} Let notations and assumptions be as before. The degree of
$V$ is defined as the maximum of the degrees of the intersections of $V$
with affine linear spaces belonging to ${\cal D}_V$. We denote the degree of
$V$ by $deg\ V$. \end{defn}

As observed in \cite{joos:tesis}, this  definition of the degree never
yields infinity, but gives always a natural number. Our definition is
equivalent to the following one: Consider all linear changes of coordinates
in $\C^n$ defined by non-singular matrices with integer entries, i.e. linear
changes of the type

$$(\xon) \mapsto (\yon),$$ 

where $Y_k = a_{k1} X_1 \plp a_{kn} X_n$ is a linear form with integer
coefficients for $1\leq k\leq n$. Any generic linear change of coordinates
induces an integral ring extension as follows:

$$\Q[Y_1,\ldots,Y_{n-i}] \longrightarrow \Qyon / (\fo i ) = \Q[V].$$

It turns out that the ring $\Q[V]$ is a free $\Q[Y_1,\ldots,Y_{n-i}]$-module
of finite rank (cf. \cite{gihesa} for instance). This rank is the same for
any generic linear coordinate change and it equals the degree of the variety
$V$. \spar



In order to define the height of an affine diophantine variety, we also
consider first the zero-dimensional case and then the case of positive
dimension.  
\spar


To start with let us first say what we mean by the (logarithmic) height of
an integer, a vector of integers, a matrix over $\Z$ and a polynomial with
integer coefficients. Let $a\in\Z$ be an integer, then the {\em height\/} of
$a$ is defined as $ht(a):= max \{log_2\abs a , 1\}$. It is obvious that the
height measures the bit length of $a$. On the other hand we shall see soon
that this simple notion of height of an integer has a natural extension to
algebraic varieties where it plays a r\^ole of ``arithmetic degree" (see
\cite{bgs2,BeYg3,Faltings,krick-pardo1,krick-pardo:CRAS,philippon1,philippon2,philippon3}).  The main outcome of this paper will be the
reinterpretation of the notion of height for algebraic varieties as a
measure for the bit complexity of an elimination procedure. In this sense,
our contribution justifies a posteriori Northcott's terminology of
``complexity" for the height of an algebraic variety \cite{SouleO}. \spar

For a vector of integer numbers $\alpha:= (a_1\klk a_n)\in\Z^n$, we define
the height $ht(\alpha)$ as the maximum of the heights of its coordinates. 
For a matrix $A\in {\cal M}_n(\Z)$ with integer entries, its height $ht(A)$
is defined as the height of $A$ as vector.  Similarly, we define the height
of a polynomial $f\in\Zxon$ as the height of the vector of its coefficients.

\begin{defn} Given a zero-dimensional diophantine algebraic variety
$V\subseteq \C^n$ and a linear form $U = \lambda_1X_1\plp \lambda_nX_n$ with
integer coefficients representing a primitive element $u$ of the ring
extension $\Q\longrightarrow \Q[V]$, we define the height of $V$ with
respect to $U$ as the maximum of the heights of the polynomials $q_u(T),
\rho_1^{(u)}X_1-v_1^{(u)}(T) \klk \rho_n^{(u)}X_n-v_n^{(u)}(T)$ (see
Definition \ref{geosolve}) and denote this height by $ht(V;U)$.
\end{defn}

In our first approach to find a suitable notion of height for algebraic
varieties we define the height of the given zero-dimensional variety $V$ as
the function $ht_V: \N \longrightarrow \N$ which associates to any natural
number $c\in\N$ the value $ht_V(c):= max \{ ht(V;u) : ht(u) \leq c \}$ if
the ring extension $\Q\longrightarrow \Q[V]$ has a primitive element of
height at most $c$ and which associates to $c$ the value $1$ if no such
primitive element exists. 
\spar

This notion of height is related to the hermitian norm and the denominator
of the points of the variety $V$.  In order to explain this relation, let us
introduce the following notations~:
\spar

Taking into account that the zero-dimensional variety $V$ is contained in
$\C^n$, we define the norm $\norm V$ as:

$$\norm V :=max\{ \norm \alpha  : \alpha\in V\}.$$
 
Furthermore, a natural number $d\in\N$ is called a {\sl
denominator of} $V$ if all elements in the set

$$d\cdot V:= \{ d\cdot \alpha : \alpha \in V\}$$

have algebraic integers as coordinates. The smallest denominator of $V$ will
be called {\em the\/} denominator of $V$ and denoted by $d_V$. \spar

With these notations we are able to state the following height estimation
for zero-dimensional varieties:

\begin{lem} 
There exists a universal constant $\kappa>0$ such that for any
zero-dimensional diophantine subvariety $V$ of $\C^n$ and any $c\in\N$ the
following inequality holds:

$$ht_V(c)\leq \delta^\kappa(c+log_2(n d_V\norm V)).$$

(This means the function $ht_v$ is of order $ht_V(c) =
\delta^{O(1)}\left(c+log_2(n d_V\norm V )\right)$).

\end{lem}

\begin{pf} Let $V\subset\C^n$ be a zero-dimensional diophantine variety of degree
$\delta\in\N$. The inequality is trivial for $ht_V(c)=1$. Therefore we may
suppose without loss of generality that the ring extension $\Q
\longrightarrow \Q[V]$ has a primitive element $u$ which is the image of a
linear form $U=\lambda_1X_1 +\cdots+ \lambda_nX_n \in \Zxon$ of height at
most $c$. We just show the inequality:

$$ 2^{ht(q_u(T))} \leq (2^{c+1}n d_V\norm V )^\delta$$

and leave the corresponding inequalities for $\rho_1^{(u)}X_1-v_1^{(u)}(T)
\klk \rho_n^{(u)}X_n-v_n^{(u)}(T) )$ to the reader. Since $d_V$ is a
denominator for $V$ and $u$ has integer coordinates, $d_V$ is also a
denominator for the set of algebraic numbers $u(V)$ contained in $\C$.
Therefore, if $\alpha_1\klk \alpha_{\delta} \in \C^n$ are the points of $V$ and if $T$
is a new indeterminate, the polynomial

$$f(T):=\prod_{i=1}^\delta\left(T-d_Vu(\alpha_i)\right)$$ 

has integer coefficients, i.e. $f(T)$ belongs to $\Z[T]$. Moreover, the
polynomial $f(d_V T)$ vanishes on the set $u(V)$.  Since $f$ and $q_u$ have
the same degree $\delta$ and since $q_u$ is the primitive minimal polynomial
of the image $u$ of $U$ in $\Q[V]$, the vanishing of $f(d_V T)$ on $u(V)$
implies that the existence of a non-zero integer $b\in\Z$ such that $f(T) =
b q_u(T)$ holds. Taking into account $\abs b \geq 1$ and that the
coefficients of $f$ are elementary symmetric functions in the points of the
set $u(V)\subset\C$ we deduce the following inequalities:

$$2^{ht(q_u)}\leq 2^{ht(f)} \leq (2 d_V \norm {u(V)})^\delta\leq
(2d_V)^\delta(n2^c)^\delta\norm{V}^\delta
\leq(2^{c+1}nd_V\norm V)^\delta.$$
\qed
\end{pf}

Let $V$ be again a zero-dimensional diophantine subvariety of $\C^n$ of
degree $\delta$. We estimate now the quantities $d_V$ and $\norm V$. 

\begin{prop}\label{ht-cota} Let $c\in\N$ be a natural number such
that the ring extension $\Q\longrightarrow\Q[V]$ has a primitive element
$u:=\lambda_1X_1+\cdots+\lambda_nX_n \in\Zxon$ of height at most $c$. Then
$d_V$ and $\norm V$ can be estimated as follows:

\begin{enumerate}

\item $$ d_V\leq \abs{ a_{\delta}^{{\delta}-1}\prod_{i=1}^n\rho_i^{(u)}},$$

where $a_{\delta}$ is the leading coefficient of $q_u(T)$ and
$\prod_{i=1}^n\rho_i^{(u)}$ is the discriminant obtained from the
polynomials in Definition \ref{geosolve}.  It follows that $V$ has
denominators of height $O((n+{\delta})ht_V(c))$.

\item the norm $\norm V$ of $V$ satisfies 

$$\norm V \leq \sqrt{n}{\delta}2^{\delta ht_V(c)}.$$

\end{enumerate}
\end{prop}

\begin{pf} The integer $\abs {a_\delta}$ is a denominator of the set $u(V)$
because it is the leading coefficient of the polynomial $q_u(T)$ which
defines $u(V)$. Since the polynomials $\rho_1^{(u)}X_1-v_1^{(u)}(T) \klk$
$\rho_n^{(u)}X_n-v_n^{(u)}(T)$ ``parametrize" the variety $V$ in function of
the zeroes of $q_u(T)$, i.e. in function of the set $u(V)$, we deduce that
$\abs{a_\delta^{\delta-1}\prod_{i=1}^n \rho_i^{(u)}}$ is a denominator of
$V$ (notice that the degree of the polynomials $v_1^{(u)}(T)\klk
v_n^{(u)}(T)$ is bounded by $\delta -1$).
\spar

>From \cite[Chapitre IV, Th\'eor\`eme 2 (ii)]{Mignotte} one deduces~:

$$\norm {u(V)} \leq 2^{ht_V(c)}.$$

Thus, for $\alpha=(\alpha_1,\ldots,\alpha_n)\in V$ and $1\leq k \leq n$ we
have:

$$\abs {\alpha_k} =\abs{\rho_k^{(u)}}^{-1} \abs {v_k^{(u)}(u(\alpha))}.$$

Since $\rho_k^{(u)}$ is an integer and $v_k(T)$ is a polynomial of degree
at most $\delta -1$ and of height at most $2^{ht_V(c)}$, we conclude~:

$$\abs{\alpha_k} \leq {\delta} 2^{ht_V(c)} \norm {u(V)}^{\delta-1} 
\leq {\delta}2^{\delta ht_V(c)}.$$

This implies

$$\norm V := max \{ \norm \alpha : \alpha\in V\} \leq
\sqrt{n}{\delta}2^{\delta ht_V(c)}.$$
\end{pf}

\begin{remark}\label{remark11}

The height is deeply related to complexity issues in polynomial equation
solving. In fact, given a regular sequence of polynomials $\fon\in\Zxon$
defining a zero-dimensional diophantine variety $V$, we have an upper bound
for the length of the output of any algorithm which solves the system $f_1=0
\klk f_n=0$ geometrically, namely:

$$(n+1)deg\ V\ ht_V(c).$$ 

Here we assume tacitly that the output is given by a linear form
representing a primitive element of the ring extension
$\Q\longrightarrow\Q[V]$ and by polynomials as in Definition \ref{geosolve} and
that these polynomials are given in dense and their coefficients in bit
representation. \end{remark}

It is also possible to give an upper bound for the height of
zero-dimensional complete intersection varieties in terms of purely
syntactical properties of its defining equations. This is the content of
what follows, namely a weak form of the so-called arithmetical B\'ezout
theorem (\cite{BeYg1,BeYg2,bgs1,bgs2,Elkadi,philippon1,philippon2,philippon3,krick-pardo1,krick-pardo:CRAS}).

\begin{prop}\label{prop-weak-bezout} Let $\fon\in\Zxon$ be
polynomials of degree at most $d$ and of height at most $h$. Suppose they
define a zero-dimensional variety $V:= V(\fon)$.  Then the function $ht_V$
can be estimated as follows: for any $c\in\N$ we have

$$ht_V(c) = d^{O(n)}\cdot h\cdot c.$$ 

\end{prop}

Proposition \ref{prop-weak-bezout} is an immediate consequence of
\cite{krick-pardo1}. Under the same assumptions as in Remark \ref{remark11},
namely that any algorithm which solves a given zero-dimensional complete
intersection polynomial equation system returns a primitive element and
polynomials as in Definition \ref{geosolve}, we can state the following
obvious {\em lower bound\/} for the bit length of the output (again we
suppose that the output polynomials as in Definition \ref{geosolve} are
given in dense and their coefficients in bit representation):

\begin{prop} 
Let $\fon\in\Zxon$ be a regular sequence defining a
zero-dimensional diophantine variety $V=V(\fon)$. Let
\linebreak
$\eta:= min \lbrace ht(V;u) : u {\hbox{ is a primitive element of the ring
extension }} \Q\longrightarrow \Q[V]\rbrace $.
Then, the bit time complexity of any algorithm 
solving $f_1=0\klk f_n=0$ geometrically is at least:
$$max \lbrace deg\ V, \eta\rbrace.$$
\end{prop}

Let us now consider the case of diophantine complete intersection variety
$V$ of positive dimension. Let $1\leq i < n$ and let be given a smooth
regular sequence of polynomials $\fo i \in \Zxon$ defining a
$(n-i)$-dimensional diophantine complete intersection variety $V$. Let
${\cal N}_V$ be the class of all linear coordinate changes:

$$ {\pmatrix { Y_1 \cr \vdots \cr Y_n } } = {\pmatrix { a_{1,1} & \cdots &
a_{1,n} \cr \vdots & \hfill{} & \vdots \cr a_{n,1} & \cdots & a_{n,n} } }
\times {\pmatrix { X_1 \cr \vdots \cr X_n } } ,$$

such that the matrix $A= (a_{ij})_{ 1\leq i,j \leq n}$ has integer entries
(i.e. $A\in {\cal M}_n(\Z)$, where ${\cal M}_n(\Z)$ denotes the ring of $n
\times n$ matrices of $\Z$), $A$ is non-singular and such that 
the ring extension

\begin{equation}\label{rex13}
\Q[Y_1,\ldots,Y_{n-i}] \longrightarrow \Qyon / (\fo i ) = \Q[V]
\end{equation}

is integral. The linear coordinate change given by the matrix $A$ induces a
{\em finite\/} morphism of affine varieties $\pi: V \longrightarrow \C^{n-i}$.
We consider for any linear form $U=\lambda_{n-i+1}X_{n-i+1}\plp
\lambda_nX_n$ with integer coefficients such that its image $u$ in $\Q[V]$ is a
primitive element of the ring extension above, the class ${\cal F}_u$ of all
those points $a\in\Z^{n-i}$ which have unramified $\pi$-fiber and for which
the image of $u$ in the coordinate ring $\Q[V_a]$ of the fiber $V_a:=
\pi^{-1}(\lbrace a\rbrace)$ remains a primitive element of the ring
extension $\Q\longrightarrow\Q[V_a]$. 
\spar

In a more down-to-earth language this means that we ask all the elements of
$V_a = \pi^{-1}(a)$ to be smooth points of $V$ and the number $deg\ V_a$ of
elements of $V_a$ to be equal to the rank of the free
$\Q[Y_1,\ldots,Y_{n-i}]$-module $\Q[V]$.  Furthermore we ask the linear form
$U$ to generate a primitive element of the ring extension
$\Q\longrightarrow\Q[V_a]$. (Observe that everything comes together since
the smoothness of the elements of $V_a$ implies that the zero-dimensional
$\Q$-algebra $\Q[V_a]$ is reduced.) Maintaining these notations we model our
notion of height of algebraic varieties by means of the following function:

\begin{defn} Given $V$ a complete intersection diophantine variety as
before, a linear change of coordinates $A\in{\cal N}_V$, a linear form
$U\in\Z[Y_{n-i+1}\klk Y_n]$ whose image $u$ in $\Q[V]$ is a primitive
element of the integral ring extension (\ref{rex13}) and given a point $a\in
{\cal F}_u$, we define the height of $V$ with respect to the triple
$(A,U,a)$ as:

$$ht(V; (A,U,a) ) := ht(V_a ; U).$$

\end{defn}

Our first approximation to the notion of height of a diophantine complete
intersection variety is given by the function $ht_V: \N \longrightarrow \N$
which associates to any natural number $c\in\N$ the value $ht_V(c):=max
\lbrace ht(V; (A,U,a)) \mid ht(A,U,a) \leq c\rbrace$ if the triple $(A,U,a)$
satisfies the assumptions of this definition, and which associates to $c$ the
value $1$ if not.

\end{subsection} 


\begin{subsection}{Straight--line Programs}\label{slp} 

In the previous subsections we discussed the mathematical form and
syntactical encoding of the output of an algorithm for geometric solving we
are going to exhibit. However we also need a suitable encoding for input and
intermediate results of our algorithm. As mathematical objects our inputs
are polynomials with integer coefficients $f\in\Zxon$. These polynomials can
be written as lists of monomials and this yields the first two possible
encodings for input and intermediate results: dense and sparse encoding. 
Thus, if $d$ is the degree of $f$ and if $h$ is an upper bound for its
height, the length of $f$ under dense encoding is

$$h \cdot {\bino{d+n}{n}} = h\cdot {\bino{d+n}{d}}.$$

Let us remark that the binomial coefficients appearing in this expression
are polynomial in the number of variables in the case of ``small" degree
polynomials (e.g. the length of the dense encoding of $f$ is $h\cdot
{\bino{n+2}{2}} \simeq h\cdot n^2$ if the degree of $f$ is at most 2).
Analogously the length of the dense encoding of $f$ is polynomial in $d$ for
$d>>n$. \spar

Alternatively, our input polynomial $f$ may be given by representing just
all its non-zero coefficients: this is the sparse encoding of the polynomial
$f$. Then, if $h$ is a bound of the height, $d$ the degree and $N$ the
number of monomials with non-zero coefficients of $f$, the length of the
sparse encoding of $f$ becomes

$$h\cdot N \cdot n\cdot log_2 d.$$

However, in many practical applications our input polynomials may be given
as programs. This is for instance the approach in the encoding of the inputs
of the version of the 3SAT problem (Problem \ref{problem1}) we stated in the
introduction. In the sequel we shall mainly restrict ourselves to
the straight--line program encoding of polynomials in the following sense:

\begin{defn} A generic straight--line program $\Gamma'$ over $\Z$ is a
pair $\Gamma' = ({\cal G}, Q)$, where ${\cal G}$ is a directed acyclic graph
and $Q$ is an assignment of instructions to the gates (i.e. vertices) of the
graph. ${\cal G}$ contains $n+1$ gates of indegree $0$ which are $Q$-labeled
by the variables \xon and by the constant $1\in\Z$. They are called the {\sl
input gates} of $\Gamma'$. We define the depth of a gate $\nu$ of the graph
${\cal G}$ as the length of the longest path joining $\nu$ and some input
gate. Let us denote the gates of the directed acyclic graph ${\cal G}$ by
pairs of integer numbers $(i,j)$, where $i$ represents the depth of the gate
and $j$ is the corresponding value of an arbitrary numbering imposed on the
set of gates of depth $i$ (this encoding can be seen in 
\cite{gihemopar,krick-pardo1,mo-mo-par,Mo-Pardo:Net,pardo}). 
Associated to the gate $(i,j)$ we have the following operation:

$$Q_{i,j} := ( \sum_{0\leq r\leq i-1} A_{i,j}^{r s}Q_{r s}) \cdot
(\sum_{0\leq r'\leq i-1} B_{i,j}^{r' s'}Q_{r' s'}),$$

where $ A_{i,j}^{r s}, B_{i,j}^{r' s'}$ are indeterminates called 
parameters of $\Gamma'$ and $Q_{rs}, Q_{r' s'}$ are precomputed values 
corresponding to the gates $(r,s)$ and $(r',s')$.

\end{defn}

We denote by $ \bar{A} = (A_{i j}^{r s}), \bar{B} = (B_{i j}^{r' s'})$ the
list of all parameters in the straight--line program $\Gamma'$.  The
intermediate results $Q_{ij}$ of $\Gamma'$ are therefore polynomials
belonging to $\Z[\bar{A},\bar{B},\xon]$ and $\Gamma'$ represents a procedure
which evaluates them. A (finite) set of polynomials $\fo s \in\Zxon$ is
said to be evaluated by a straight--line program of generic type $\Gamma'$
with parameters in a set ${\cal F}\subseteq \Z$ if specializing the
coordinates of the parameters $\bar{A}$ and $\bar{B}$ in $\Gamma'$ to values
in ${\cal F}$, there exist gates $(i_1,j_1)\klk (i_s,j_s)$ of $\Gamma'$ such
that

$$f_k = Q_{i_k j_k}(\bar{a}, \bar{b}, \xon)$$

holds for $1\leq k\leq s$.  Specializing in the indicated way the parameters
of $\Gamma'$ into values of ${\cal F}$ we obtain a copy $\Gamma$ of the
directed acyclic graph ${\cal G}$ underlying the generic straight--line
program $\Gamma'$ and of its instruction assignment $Q$. We call this copy
$\Gamma$ a {\em straight--line program (of generic type $\Gamma'$) in
$\Zxon$ with parameters in ${\cal F}$}. The gates of $\Gamma$ correspond to
polynomials belonging to $\Zxon$. These polynomials are obtained from the
intermediate results $Q_{ij}$ of $\Gamma'$ by specializing adequately the
parameters from which the depend. We shall call these polynomials the
intermediate results of the straight--line program $\Gamma$. Furthermore the
polynomials $\fo s$ are called {\em (final) results\/} or {\em outputs\/} of $\Gamma$. 
Alternatively, we shall say that $\fo s$ are represented, computed or
evaluated by $\Gamma$.\spar

The current complexity measures for the generic straight--line
program $\Gamma'$ are:

\begin{itemize} 
\item the size of $\Gamma'$ = the size of the graph ${\cal G}$ 
\item the non-scalar depth ($\Gamma'$) = the depth of the graph ${\cal G}$ 
\end{itemize}

The size and non-scalar depth of the specialized straight--line program $\Gamma$
are analogously defined. Additionally, for any straight--line program
$\Gamma$ in \linebreak
$\Zxon$ with parameters in ${\cal F}\subseteq \Z$, we call the the
maximum of the heights of the elements in ${\cal F}$ {\em the height of the
parameters of $\Gamma$\/} or, for short, the {\em height of $\Gamma$}.
The encoding of a polynomial $f\in\Zxon$ by a straight--line program
$\Gamma$ in $\Zxon$ with parameters of height $h$ has (bit) length

$$2 (L^2 h + L log_2 L).$$

Note that the notion of straight--line program encoding covers as well both
notions of dense and sparse encoding of polynomials. This can be seen as
follows: a polynomial of degree at most $d$ and height at most $h$, given in
dense encoding, can be evaluated by a straight--line program of size $O(d
{\bino{d+n}{n}})$ and non-scalar depth $O(log_2 d)$ with parameters of
height at most $h$. Similarly, if a polynomial $f\in\Zxon$ has degree at
most $d$, height $h$ and $N$ non-zero coefficients, it can be evaluated by a
straight--line program of size $O(d\cdot N)$ and non-scalar depth $O(log_2
d)$ with parameters of height at most $h$.

\end{subsection} 


\begin{subsection}{Complexity of geometric solving}\label{complejidad}

Now that the notions of geometric solving, degree and height of
polynomial equation systems and straight--line programs and its distinct
complexity measures have been introduced, we are able to state our main
result, namely the following Theorem \ref{maquina}. This theorem represents
our principal contribution to the solution of Problem \ref{problem3} in
Section \ref{intro}.

\begin{thm}\label{maquina} There exists a bounded error probabilistic
Turing machine that from a smooth regular sequence $\fon\in\Zxon$ outputs a
linear form $U\in\Zxon$ representing a primitive element $u$ of the ring
extension $\Q\longrightarrow \Q[V]$, where $V:=V\ifon$ is the algebraic
variety defined by $\fon$, and polynomials of the form $q_u(T),
\rho_1^{(u)}X_1-v_1^{(u)}(T)\klk \rho_n^{(u)}X_n-v_n^{(u)}(T)$ with
$\rho_1^{(u)}\klk \rho_n^{(u)}$ non-zero integers and $q_u, v_1^{(u)}\klk
v_n^{(u)}\in\Z[T]$ such that these polynomials represent a geometric
solution of the equation system $f_1=0\klk f_n=0$ (see Definition
\ref{geosolve}). The Turing machine finds this solution in time (counting
the number of bit operations executed)

$$(n d h L \delta \eta)^{O(log_2 n +\ell)}$$

using only 

$$(n d L \delta)^{O(1)}$$

arithmetic operations in $\Z$ at unit cost. We assume that the
polynomials \fon have degree at most $d$ and that they are given by a
straight--line program of size $L$ and non-scalar depth $\ell$ with
parameters of height at most $h$. The quantity $\delta$ defined by $\delta:=
max\{deg\ V(\fo i ): 1\leq i\leq n\}$ is the degree of the system $f_1=0\klk
f_n=0$. Finally the quantity $\eta$ is the height of the system $f_1=0\klk
f_n=0$ defined as $\eta:= max \{ ht_{V_i}(c_3 ( (log_2 n +\ell) \cdot
log_2\delta) ): 1\leq i\leq n\}$ where $c_3>0$ is a suitable universal
constant independent of the specific input $\fon$ (or even its size).
\end{thm}

Theorem \ref{maquina} follows from the description of the algorithm given in
Section \ref{algoritmo} below.  Theorem \ref{teouno} above follows
immediately from Theorem \ref{maquina} above by putting $d:=2$, $L:=n$,
$\ell:=2$. Let us remark here that Theorem \ref{maquina} improves and
extends the main result of \cite{gihemomorpar} and \cite{gihemopar} to the
bit complexity model. It sheds also new light on the main complexity outcome
of the papers \cite{ShSm1,Emiris,CaEm95}.

\end{subsection} 


\begin{subsection}{A Division Step in the Nullstellensatz}\label{step}

In this subsection we show how Nullstellensatz bounds imply Liouville
estimates. This establishes a close connection between Problems
\ref{problem4} and \ref{problem5} in the Introduction.  Let us recall the
assumptions in the statement of Theorem \ref{Liouville}. \spar

There is given a smooth regular sequence of polynomials $\fon\in$\linebreak
$\Zxon$ of degree at most $d$ and of height at most $h$, defining a
zero-dimensional affine variety $V:= V\ifon\subseteq \C^n$. Moreover there
is given a point $\alpha=(\alpha_1\klk \alpha_n)\in V$ and a real number
$0<\varepsilon\leq 1$. Finally, we assume that $V\cap \Q[i]^n = \emptyset$
holds. By a well-chosen linear change of coordinates we may assume without
loss of generality $V\cap (\Q[i]\times \C^{n-1})=\emptyset$.  It is
sufficient to show that the algebraic number $\alpha_1$ is hard to
approximate in the sense of the conclusion of Theorem \ref{Liouville}. For
this purpose consider integers $p\in \Z[i]$ and $q\in \N$ such that $p\over
q$ is an approximation at level $\varepsilon$ to the algebraic number
$\alpha_1$. \spar

We introduce the following polynomial~:

$$ f_{n+1}:= (qX_1-p)(qX_1-\bar{p})$$

(here $\bar{p}$ stands for the complex conjugate of $p$). The assumption
$V\cap (\Q[i]\times \C^{n-1})=\emptyset$ implies that the sequence of
polynomials $f_1,\ldots, f_{n+1}$ has no common zero in $\C^n$.  Therefore
there exists a non-zero integer $b\in\Z$ and polynomials $g_1,\ldots,g_{n+1}
\in \Z[X_1,\ldots,X_n]$ such that the following B\'ezout identity holds~:

$$b=g_1f_1+\cdots + g_{n+1}f_{n+1}.$$

Evaluating this identity at the point $\alpha\in V$ we obtain~:

$$1\leq \abs b = \abs{ g_{n+1}(\alpha)} \cdot\abs {q\alpha_1 -p} \cdot
\abs {q \alpha_1- \bar{p}}.$$

This yields the following estimate:

\begin{prop}\label{estrategia}

With assumptions and notations as before, for any rational approximation
$a:= (p_1/q,\ldots,p_n/q)\in\Q[i]^n$ at level $\varepsilon$ to the point
$\alpha \in V$ we have the inequality~:
 
$$ { {1} \over { \abs {g_{n+1}(\alpha)} q^2}} \leq
\norm {a- \alpha} ( \abs {\alpha_1} + \abs {{{p_1} \over {q}}} ).$$

In particular, for $\norm {a- \alpha} < \varepsilon \leq 1$, we
obtain the estimation~:

$$ { {\varepsilon^{-1}}\over {\abs {g_{n+1}(\alpha)} (2\norm {\alpha
}+1)}}\leq q^2.$$

\end{prop} 

{\em (Here $p_1\klk p_n\in\Z[i]$ are Gaussian integers and $q\in\N$ is a
natural number.)}\spar

>From the second inequality we deduce that it is sufficient to bound the
values of $\norm \alpha $ and $\abs {g_{n+1}(\alpha)}$ in order to obtain a
Liouville estimate for the rational approximation $a$ of the point $\alpha$. 
Note that the bounds for $\norm \alpha $ and $\abs {g_{n+1}(\alpha)}$ which
can be easily deduced from the known Nullstellens\"atze (as e.g. in
\cite{BeYg1,BeYg2,krick-pardo1}) would be insufficient for
our purpose (proving Theorem \ref{Liouville}) as they imply only a
unspecific general estimate, namely

$$ { {-log_2 \varepsilon} \over {d^{Cn}h } } < log_2 \abs q $$ 

for a suitable constant $C>0$. A Liouville estimate of this type does not
take into account the specific properties of the variety $V$ expressed
through its degree and height.  The more specific bounds for $\norm \alpha $
and $\abs {g_{n+1}(\alpha)}$ as required for the proof of Theorem
\ref{Liouville} by means of Proposition \ref{estrategia} are an immediate
consequence of the following result:

\begin{thm}\label{paso-division}
Let $\fon, f_{n+1}\in\Zxon$ polynomials having no common zero in $\C^n$
and verifying the following assumptions:

\begin{itemize}
\item there exists a straight--line program of size $L$ and non-scalar depth $\ell$
      with parameters of height $h$ that evaluates the polynomials
      $\fon, f_{n+1}$ 

\item the degrees of the polynomials $\fo {n+1}$ are bounded by $d$ and 
      $h'$ is an upper bound for their heights
\end{itemize}

Let us furthermore assume that $\fon$ form a smooth regular sequence which
defines a zero-dimensional affine algebraic variety $V=V\ifon$. We also
consider the following quantities:

\begin{itemize}
\item $\delta:=deg(V)$ 
\item $\eta:=$ min ht $V$ (i.e. $\eta$ is the minimal value of
$ht_V$ distinct from $1$)
\end{itemize}

Then there exists a straight--line program of size $L (nd\delta)^{O(1)}$ and
non-scalar depth of order $O(log_2 n + log_2 d + log_2\delta + \ell)$ with parameters
of height at most $O(max\lbrace h, h', \eta, log_2 n, \ell\rbrace)$ which
evaluates a non-zero integer $a\in\Z$ and a polynomial $g_{n+1}\in\Zxon$
such that

$$a-g_{n+1}\cdot f_{n+1}$$

belongs to the ideal \ifon generated by \fon in \Zxon.
\end{thm}

The proof of this theorem will follow from the description of the algorithm
given in Section \ref{division} below.  Applying Theorem \ref{paso-division}
and Lemma \ref{cotas} below, we obtain a more precise upper bound for the
value $\abs {g_{n+1}(\alpha)}$. From this bound together with Proposition
\ref{ht-cota} and \ref{estrategia} we then deduce easily Theorem
\ref{Liouville}.

\end{subsection} 

\end{section} 


\typeout{Section 3}

\begin{section}{An algorithm for Geometric Solving}\label{algoritmo} 

The aim of this section is to establish a proof for Theorem \ref{maquina}.
We describe an algorithm which implies Theorem \ref{maquina}. This algorithm
works inductively on the codimension of the varieties $V_i:=V\ifo i$, $1\leq
i\leq n$, and our main goal is to describe this recursion.  \spar

Recall that our input is a smooth regular sequence 
$f_1,\ldots,f_n \in \Z[X_1,\ldots,X_n]$ of degree at most $d$. We assume
that this input is encoded by a straight--line program $\Gamma$ of size $L$
and non--scalar depth $\ell$ with parameters of height at most $h$ that
evaluates the polynomials $f_1,\ldots,f_n$. Our algorithm computes a
geometric solution of the zero-dimensional algebraic variety $V:=
V(f_1,\ldots,f_n) \subseteq \C^n$. In order to describe for $1\leq i\leq n$
the $i$-th recursive step of our algorithm, we shall refer to the
intermediate complete intersection algebraic varieties as
$V_i:=V(f_1,\ldots,f_i)$ and introduce the following parameters:

\begin{itemize} 
\item $\delta_i:=deg(V_i)$, 
\item $\delta:= max\{ \delta_i\; : \; 1\leq i \leq n \}$,
\item $\eta_i:=ht_{V_i}(C)$, where $C$ is a suitably chosen natural number
of order \linebreak
$O((log_2n+\ell)log_2\delta)$, such that $\Q[V_i]$ has a primitive
element of height $C$ with respect to a suitable Noether position of $V_i$.

\item $\eta:=max\{\eta_i\;:\; 1\leq i \leq n\}$. 
\end{itemize}

As a byproduct of our algorithm, we shall obtain geometric solutions of the
algebraic varieties $V_i$. In fact, our algorithm computes for every 
$1\leq i \leq n$  a $\Z$-linear change of coordinates

$$(X_1,\ldots,X_n) \longrightarrow (Y_1^{(i)},\ldots,Y_n^{(i)}),$$

such that the ring extension

$$R_i:=\Q[Y_1^{(i)},\ldots,Y_{n-i}^{(i)}]\longrightarrow
B_i:=\Q[Y_1^{(i)},\ldots,Y_n^{(i)}]/(f_1,\ldots,f_i)= \Q[V_i]$$

is integral.  The $R_i$-module $B_i$ is a free module of rank say $D_i\leq
\delta_i$. Let us denote by $\pi_i:V_i\longrightarrow \C^{n-i}$ the
projection on the first $n-i$ coordinates of $(Y_1^{(i)},\ldots,Y_n^{(i)})$.
The morphism $\pi_i$ of affine varieties is finite.

\begin{defn} Let assumptions and notations be as before.  A {\sl
lifting point for $W:=V_i$} of the finite morphism $\pi_i$ is a point $P =
(p_1\klk p_{n-i})\in\Z^{n-i}$ with the following properties:

\begin{itemize} 

\item the zero-dimensional fiber $W_{P}:=\pi_i^{-1}(P)$ has degree (i.e.
cardinality) equal to the rank of $B_i$ as free $R_i$-module (this means
$deg(W_{P})=D_i$)

\item the fiber $W_{P}$ contains only smooth points. (This is equivalent to
saying that the Jacobian matrix of the polynomials we obtain from
\linebreak 
$f_1(Y_1^{(i)}\klk Y_n^{(i)}) \klk f_i(Y_1^{(i)}\klk Y_n^{(i)})$ by
substituting for $Y_{1}^{(i)}\klk Y_{n-i}^{(i)}$ the coordinates $p_1\klk
p_{n-i}$ of $P$ is regular at every point of the fiber $W_{P}$.)

\end{itemize}

If $P\in\Z^{n-i}$ is a lifting point of the morphism $\pi_i$ we will call
its fiber $W_{P}$ a {\em lifting fiber\/} of $W=V_i$. Observe that the
elements of a lifting fiber of $\pi_i$ are smooth points of $W$. 
\end{defn}

The lifting fibers $W_{P}$ have the property that a geometric solution of
the variety $V_i$ can be reconstructed from the projection $\pi_i$ and any
geometric solution of the equations of such a fiber, (see Subsection \ref{lifting} below).
Our algorithm will choose for each $1\leq i\leq n$ a suitable lifting point
$P_i\in\Z^{n-i}$ of the morphism $\pi_i$. In the sequel we shall denote the
lifting fiber of this point by $V_{P_i}$.
\spar

The following Section \ref{algoritmo} is divided into two two
well--distinguished parts, namely: 

\begin{enumerate} 

\item Subsection \ref{lifting}, treating ``lifting by a symbolic Newton
method", where we show how it is possible to reconstruct a geometric
solution of the equations $\fo i$ from the lifting fiber $V_{P_i}$.

\item Subsection \ref{liftingpoint}, where  we show how to find a linear
coordinate change 

$(\xon) \longrightarrow (Y_1^{(i)} \klk Y_n^{(i)})$, the
lifting point $P_i$ and a geometric solution of the equations of the lifting
fiber $V_{P_i}$.

\end{enumerate}

The $i$-th recursive step of our algorithm consists of the computation of
the new variables $Y_1^{(i+1)},\ldots,Y_n^{(i+1)}$, the point $P_{i+1}$ and
a geometric solution of its fiber $V_{P_{i+1}}$ from the data
$Y_1^{(i)},\ldots,Y_n^{(i)}$, $P_i$ and $V_{P_i}$. \spar

Let us conclude this subsection with the following remark: let $1\leq i\leq
n$ and consider the lifting fiber $V_{P_i}$. Since the degree of $V_{P_i}$
equals $D_i= rank_{R_i}(B_i)$, any $\Z$-linear form
$U_i:=\lambda_{n-i+1}Y_{n-i+1}^{(i)}+\cdots +\lambda_{n}Y_n^{(i)}$ which
separates the points of $V_{P_i}$ represents not only a primitive element of
the ring extension $\Q\longrightarrow \Q[V_{P_i}]$ but also a primitive
element of the integral ring extension $R_i\subseteq B_i$ (which we will
denote by $u_i$). Our algorithm will compute a geometric solution of both
$V_i$ and $V_{P_i}$ using such a linear form $U_i$. \spar

Let us observe that the bit length of our geometric solution of the
equations of $V_{P_i}$ will be bounded by the quantity

$$(i+1)(\delta_i+2)\eta_i.$$


\begin{subsection}{Elementary operations and bounds for straight--line
programs}\label{elementary}

In this subsection we collect some elementary facts about straight--line
programs. We start with an estimate for the degree and height of a
polynomial given by a straight--line program. In order to state our result
with sufficient generality, let us observe that the notion of height makes
sense mutatis mutandi for polynomials over any domain equipped with an
absolute value.

\begin{lem}[\cite{krick-pardo1}]\label{cotas} Let $R$ be a ring equipped
with an absolute value $\abs \cdot : R \longrightarrow \R$. Suppose
$f\in\Rxon$ is a polynomial which can be evaluated by a straight--line
program $\Gamma$ in $\Rxon$ of size $L$ and non-scalar depth $\ell$
with parameters of height $h$. Let $H>0$ a real number and let 
$\alpha = (\alpha_1,\ldots,\alpha_n)$ be a point in
$R^n$ such that $log_2 \abs {\alpha_i} \leq H$ holds for $1\leq i\leq n$. 
Then $f$ and $f(\alpha)$ satisfy the following estimations:

\begin{itemize} 
\item $deg(f)\leq 2^{\ell}$ 
\item $ht(f)\leq (2^{\ell+1}-1)\cdot (h + log_2 L)$ 
\item $log_2 \abs{f(\alpha)} \leq (2^{\ell+1}-1)\cdot (max \{ h,H \} +  log_2 L)$ \end{itemize}
\end{lem}

One of the main ingredients used in our procedure below is the efficient
computation of the coefficients of the characteristic polynomial of a
matrix. We shall use for this task Berkowitz' division-free and well
parallelizable algorithm \cite{berk} (compare also \cite{Csanky} and
\cite{BoGaHo} for historical predecessors of this algorithm). This is the
content of the the next lemma.

\begin{lem}\label{berk}(\cite{berk,krick-pardo1}) Let $R$ be a
domain. There exists a straight--line program of size $N^{O(1)}$ and
non-scalar depth $O(log_2 N)$ with parameters in $\lbrace -1, 0, 1\rbrace$
that from the entries of any $N\times N$ input matrix over $R$ computes all
coefficients of the characteristic polynomial of the given matrix.
\end{lem}

We use this algorithm not only for the computation of the characteristic
polynomial of a given matrix but also for the computation of the greatest
common divisor of two given univariate polynomials with coefficients in a
unique factorization domain (this task can be reduced to solving a suitable
linear equation system corresponding to the B\'ezout identity over the
ground domain. See \cite{krick-pardo1} for details).\spar

The following result is an immediate consequence of the formal rules of
derivation:

\begin{lem}\label{diff} Let $R$ be a domain. For a given a finite set of
polynomials $f_1 \klk f_s$ of $\Rxon$ which can be evaluated by a
straight--line program $\beta$ in\linebreak 
$\Rxon$ of size $L$ and non-scalar depth
$\ell$, there exists a straight--line program in $\Rxon$ of size $(2n+1)L$
and non-scalar depth $\ell +1$ with the same parameters as $\beta$ which
evaluates $ f_1\klk f_s$ and all the first partial derivatives:

$$ \lbrace {{\partial f_i} \over {\partial X_j}} : {1\leq i\leq s}, {1\leq
j\leq n} \rbrace.$$

\end{lem}

Combining this lemma with Lemma \ref{berk}, one concludes: let $f_1\klk f_n$
be a family of polynomials of $\Rxon$ which can be evaluated by a
straight--line program $\beta$ of size $L$ and non-scalar depth $\ell$. Then
there exists a straight--line program in $\Rxon$ of size $n^{O(1)}L$ and
non-scalar depth $O(\ell+log_2n)$ with the same parameters as $\beta$ which
evaluates the Jacobian determinant.

$$J\ifon : = det ({{\partial f_i} \over {\partial X_j}})_{1\leq i,j\leq n}.$$ 

In some exceptional cases the straight--line programs we are going to
consider might contain divisions as operations. Since we are only interested
in division-free straight--line programs, the following ``Vermeidung von
Divisionen" technique due to V. Strassen \cite{stras:verm} becomes crucial:

\begin{prop}[\cite{krick-pardo1}, \cite{stras:verm}]\label{verm} Let
$\Gamma$ be a (division-free) straight--line program in $\Zxon$ of size $L$
and depth $\ell$ with parameters of height $h$ that computes a finite set of
polynomials $f_0,\ldots,f_m$ of $\Z[X_1,\ldots, X_n]$. Assume that $f_0\neq
0$ holds and that $f_0$ divides $f_i$ in $\Z[X_1,\dots, X_n]$ for any $1\leq
i \leq m$. Then there exists a (again division-free) straight--line program
$\Gamma'$ in $\Zxon$ with the following properties:

\begin{enumerate}

\item  $\Gamma'$ computes  polynomials $P_1,\ldots, P_m$ of
$\Z[X_1,\ldots,X_n]$ and a non-zero integer $\theta$ such that
for any $1\leq i\leq m$ holds

$$ P_i=\theta {{f_i}\over {f_0}}.$$

\item $\Gamma'$ has size of order $  O(2^{2\ell}(L+n+2^{\ell}+m))$,
depth of order $O(\ell)$  and its parameters have height of order 
$\max\{h,O(\ell)\}$.
\end{enumerate}

Moreover, the height of $\theta$ is of order:

$$2^{O(\ell)}( max\{ h, \ell\} + log_2 L).$$

\end{prop}

The proof of this proposition is based on the computation of the homogeneous
components of a polynomial given by a straight--line program. This proof
also provides an algorithm computing the homogenization of a polynomial
given by a straight--line program. This is the content of the next lemma.

\begin{lem}[\cite{krick-pardo1}]\label{comphom} Suppose that we are given
a polynomial\linebreak
$P:=\sum P^\mu\, X_1^{\mu_1}\dots X_n^{\mu_n}$ in
$\Z[X_1,\dots,X_n]$ which can be evaluated by a straight line program
$\Gamma$ of size $L$ and depth $\ell$ with parameters in a given set ${\cal
F}\subset \Z$. Let be given a natural number $D$. Then there exists a
straight line program $\Gamma '$ in $\Zxon$ which computes all the
homogeneous components of $P$ having the following properties: $\Gamma'$
uses parameters from ${\cal F}$, has size $(D+1)^2$ and non-scalar depth
$2\ell$.  \end{lem}

In Section \ref{division} we shall work with a specific polynomial which we
call the pseudo-jacobian determinant of a given regular sequence.  We
introduce now this polynomial and say how it can be evaluated. Let $R$ be a
domain containing $\Q$. Let $\K$ be the field of fractions
of $R$ and let $\fon\in \Rxon$ be a regular sequence in $\Kxon$. 
Furthermore, let $Y_1\klk Y_n$ be new variables. We write $Y=(Y_1\klk Y_n)$.
Fix $1\leq j\leq n$.  By $f_j^{(Y)} := f_j(Y_1\klk Y_n)$ we denote the
polynomial from $f_j$ substituting the variables $\xon$ by $\yon$. In the
polynomial ring $R[\yon, \xon]$ we decompose the polynomial $f_j^{(Y)}-f_j$
in the following (non-unique) way: 

$$ f_j^{(Y)}-f_j = \sum_{k=1}^n l_{k,j}(Y_k-X_k),$$ 

with $l_{k,j} \in R[Y_1\klk Y_n,\xon]$. Let us consider the determinant
$\Delta$ of the matrix $A = (l_{k,j})_{1\leq k,j\leq n}$, namely:

$$\Delta:= det(A).$$

This determinant is called the {\em pseudo-jacobian determinant\/} of the
regular sequence of polynomials $\fon$. If $d$ is a bound for the degrees of
\fon and these polynomials are given by a straight--line program $\beta$ of
size $L$ and non-scalar depth $\ell$, then there is a straight--line program
$\beta'$ of size $(nd)^{O(1)}L$ and non-scalar depth $O(log_2 n + \ell)$
which evaluates the pseudo-jacobian determinant $\Delta$.  The
straight--line program $\beta'$ uses apart from the same parameters as
$\beta$ only parameters of $\Z$ of height $O(log_2d)$.
\spar

We shall also consider the execution of straight--line programs in matrix
rings. The situations where we apply these considerations will be of the
following type: Let $R$ be a domain. Suppose that there is given a
polynomial $g\in\Rxon$ by a straight--line program $\beta$ in $\Rxon$ of
size $L$ and non-scalar depth $\ell$.  Suppose also that there are given $n$
commuting $D\times D$ matrices $M_1\klk M_n$ over $R$. In such a situation
the entries of the matrix $g(M_1\klk M_n)$ can be computed from the entries
of of $M_1\klk M_n$ and the parameters of $\beta$ by a straight--line
program $\beta'$ in $R$ of size $D^{O(1)}L$ and non-scalar depth $O(\ell)$.
The parameters of $\beta'$ are just the values $0,1$ (see \cite{gihesa} for
details). \spar

Another important aspect of our main algorithm is its probabilistic (or
alternatively its non-uniform) character. This is the content of the next
definition and proposition.

\begin{defn} Let be given  a set of polynomials ${\cal W}
\subseteq\Zxon$. A finite set $Q\subseteq\Z^n$ is called a correct test
sequence (or questor set) for ${\cal W}$ if for any polynomial $f$
belonging to ${\cal W}$ the following implication holds:

$$f(x) = 0\hbox{ for all } x\in Q \hbox{ implies } f = 0.$$

\end{defn}

Denote by ${\cal W}(n,L,\ell)$ the class of all polynomials of $\Zxon$ which
can be evaluated by straight--line programs of size at most $L$ and of
non-scalar depth at most $\ell$. The following result says that for the
class ${\cal W}(n,L,\ell)$ exist many correct test sequences of moderate
length.

\begin{prop} [\cite{heschnorr,krick-pardo1}]\label{questores}
Let be given natural numbers $n,\ell, L$ with $L\geq n+1$ and consider the
following quantities:

$$u:=(2^{\ell +1}-2)\, (2^{\ell} +1)^2 \qquad {\rm and} \qquad t:= 6 \,(\ell
L)^2.$$

Then the finite set $\{1,\dots,u\}^{nt}\subset \Z^{nt}$ contains at least
$u^{nt}\,(1-u^{-{ t\over 6}})$ correct test sequences of length $t$ for
${\cal W}(n,L,\ell)$.  In particular the set of correct test sequences for
${\cal W}(n,L,\ell)$ of length $t$ containing only test points from
$\{1,\dots,u\}^{n}$ is not empty.

\end{prop}

>From Lemma \ref{cotas} we deduce the following complexity estimate:
\spar

Let $f\in\Zxon$ be a polynomial given by a straight--line program of size
$L$ and non-scalar depth $\ell$ with parameters of height $h$. Let
$\alpha\in\Z^n$ be a point of height $h'$ given in bit representation. Then
there exists a (deterministic) ordinary Turing machine which computes the
bit representation of the value $f(\alpha)$ in time $(2^{\ell} L max\lbrace
h, h'\rbrace )^{O(1)}$. \spar

In the next subsection we shall make use of a problem adapted version of the
Hensel-Newton iteration.  We are now going to describe a suitable
division-free symbolic form of this procedure. \spar

Let $R$ be a polynomial ring over $\Q,$ let $K$ be its field of fractions
and let $f_1,\ldots,f_n\in R[X_1,\ldots,X_n]$ be polynomials of degree at
most $d$. Suppose that $f_1,\ldots,f_n$ are given by a (division-free)
straight--line program $\beta$ of size $L$ and non--scalar depth $\ell$. 
Let us also assume that the Jacobian matrix $D(f):=D(f_1,\ldots,f_n):=\left(
{{\partial f_i} \over {\partial X_j}} \right) _{ 1 \le i,j \le n}$ of the
polynomials $f_1,\ldots,f_n$ is regular. We consider now the following
Newton--Hensel operator:

\begin{equation}\label{uno} N_f(X_1,\ldots,X_n):=\pmatrix{X_1\cr \vdots\cr
X_n \cr} - D(f)^{-1} \pmatrix{f_1(X_1,\ldots,X_n)\cr \vdots\cr
f_n(X_1,\ldots,X_n) \cr}.\end{equation}

This operator is given as a vector of $n$ rational functions of
$K(X_1,\ldots,X_n)$. This is also true for the $k$-th iteration of this
operator, which we denote by $N_f^k$.  For any $k\in\N$ there exist
numerators $(g_1^{(k)},\ldots,g_n^{(k)}) \in R[X_1,\ldots,X_n]$ and a
non-zero denominator $h^{(k)}\in R[X_1,\ldots,X_n]$ such that $N_f^k$ can be
written as:

$$N_f^k=\left({{g_1^{(k)}}\over{h^{(k)}}},\ldots,{{g_n^{(k)}}\over{h^{(k)}}}
\right)\in K(X_1,\ldots,X_n)^n.$$

The next lemma gives a description of a division-free straight--line
program in $\Rxon$ that evaluates numerators and denominators for $N_f^k$.

\begin{lem}\label{lemanewton}  Let notations and assumptions be as before.
Let $k$ be a natural number. There exists a straight--line program in
$\Rxon$ of size $O(kd^2 n^7L)$ and non-scalar depth $O((log_2n +\ell)k)$
with the same parameters as $\beta$ which evaluates numerators
$g_1^{(k)},\ldots,g_n^{(k)}$ and a (non-zero) denominator $h^{(k)}$ for the
$k$-fold iteration $N_f^k$ of the Newton-Hensel operator $N_f$. \end{lem}

\begin{pf} Let $A(f) = (a_{ij})_{1\leq i,j\leq n}$ be the transposed matrix of the
adjoint matrix of $D(f)$. Both $A(f)$ and the Jacobian determinant
$J(f):=det(D(f))$ can be evaluated by a straight--line program of size
$O(n^5+nL)$ and non--scalar depth $O(log_2n+\ell)$, as it can be seen just
by combining Lemma \ref{berk} and \ref{diff}.  We can
write the operator $N_f$ as

\begin{equation} \label{dos} N_f={{ J(f) \pmatrix{X_1\cr \vdots\cr X_n \cr} -
A(f) \pmatrix{f_1(X_1,\ldots,X_n)\cr \vdots\cr f_n(X_1,\ldots,X_n) \cr}} \over
{J(f)}}.\end{equation}

The entries $a_{ij}$ of the matrix $A(f)$ are polynomials
of the ring $R[X_1,\ldots,X_n]$ having degree at most $(n-1)(d-1)$. 
Moreover, the Jacobian determinant $J(f)$ is a polynomial of
$R[X_1,\ldots,X_n]$ having degree at most $n(d-1)$. 
For $1\leq i\leq n$ we consider 

$$g_i:= J(f) X_i - \sum \limits_{j=1}^n a_{i,j} f_j.$$ 

All polynomials appearing on the right hand side of the definition of $g_i$
as summands have degree bounded by $\nu:=nd +1$. Thus the degree of any
$g_i$ is bounded by $\nu$. Let $^hg_i(X_0, X_1\klk X_n)\in
R[X_0,\ldots,X_n]$ be the homogenization of $g_i$ by a new variable $X_0$
and let $^hJ(f)(X_0,X_1,\ldots,X_n)\in R[X_0,\ldots,X_n]$ be the
homogenization of the Jacobian determinant $J(f)$ by $X_0$. \spar

We introduce now the following homogeneous polynomials (forms): 

\begin{itemize} 
\item $G_i(X_0,\ldots,X_n):=X_0^{\nu-deg(g_i)}({^hg_i})$, 

\item $H(X_0,\ldots,X_n):=X_0^{\nu-deg(J_f)}({^hJ(f)})$. 
\end{itemize} 

According to Lemma \ref{comphom}, there exists a division-free
straight--line program in $\Rxon$ of size $O(d^2(n^7+n^3L))$ and non--scalar
depth $O(log_2n+\ell)$ which evaluates the forms $G_1,\ldots,G_n,H$. We now
define recursively the following polynomials:

\begin{itemize}

\item{} for $k=1$, $1\leq i \leq n$ let $g_i^{(1)}:=G_i(1,X_1,\ldots,X_n)$,
$h^{(1)}:=H(1,X_1,\ldots,X_n)$,

\item{} for $k\ge 2$, $1\leq i \leq n$ let $g_i^{(k)}:=G_i(h^{(k-1)},
g_1^{(k-1)},\ldots,g_n^{(k-1)} )$,
\linebreak $h^{(k)}:=H(h^{(k-1)}, g_1^{(k-1)},\ldots,g_n^{(k-1)} ).$

\end{itemize}

It is easy to see that these polynomials $g_1^{(k)},\ldots,g_n^{(k)}$ are
numerators and that the polynomial $h^{(k)}$ is a denominator of the iterated
Newton-Hensel operator $N_f^k$. A straight--line program evaluating them is
obtained by iterating $k$ times the straight--line program which computes
$G_1,\ldots,G_n$ and $H$. No new parameters are introduced by this
procedure. Putting all this together we get the complexity bounds in the
statement of Lemma \ref{lemanewton}. \qed
\end{pf}

\end{subsection} 


\subsection{Lifting fibers by the symbolic Newton-Hensel algorithm}\label{lifting}

The idea of using a symbolic adaptation of Newton--Hensel iteration for
lifting fibers was introduced in \cite{gihemomorpar} and \cite{gihemopar}. 
For technical reasons, in these papers it was necessary to use algebraic
parameters for the lifting process. We present here a new version of this
lifting algorithm in which the use of algebraic numbers is replaced by a
certain matrix with integer entries. The whole procedure therefore becomes
completely rational. The new lifting process is described in the statement
of the next theorem and its proof. \spar

Let notations and assumptions be as the same as at the beginning of this
section. We fix $1\le i \le n$ and assume for the sake of notational
simplicity that the variables \xon are already in Noether position with
respect to the variety $V_i$, the variables $\xo {n-i}$ being free.  We
suppose that the lifting point $P_i$, the coordinates of the $\Z$-linear
form $U_i$ and a geometric solution for the equations of the lifting fiber
$V_{P_i}$ are explicitly given. With these conventions we state the main
result of this subsection as follows:

\begin{thm} \label{newton} There exists a (division-free) straight--line
program $\Gamma_i$ in the polynomial ring $\Z[\xo {n-i}, U_i]$ of size
$(id\delta_iL)^{O(1)}$ and non-scalar depth $O((log_2 i+\ell)log_2
\delta_i)$ using as parameters

\begin{itemize}
\item the coordinates of $P_i$
\item the integers appearing in the geometric solution of the equations of
the lifting fiber $V_{P_i}$ and 
\item the parameters of the input program $\Gamma$
\end{itemize} 

such that the straight--line program $\Gamma_i$ computes
\begin{itemize}

\item the  minimal polynomial $q_i \in \Z[X_1,\ldots,X_{n-i},U_i]$ of the
primitive element $u_i$ of the ring extension $\Qxo {n-i} \longrightarrow
\Q[V_i] = \Qxon /\ifo i$,

\item polynomials $\rho^{(i)}_{n-i+1},\ldots,\rho^{(i)}_n \in
\Z[X_1,\ldots,X_{n-i}]$, $\rho:= \prod\limits_{k=n-i+1}^n\rho^{(i)}_k$ and
polynomials $v^{(i)}_{r+1},\ldots,v_n^{(i)} \in \Z[X_1,\ldots,X_{n-i},U_i]$
with $max\{deg_{U_i} v_j^{(i)}\; ; \; r< j \le n\} < \delta_i$ such that

\end{itemize}

$$\ifo {i}_{\rho}=(q_i(U_i),\rho^{(i)}_{r+1}X_{r+1}-
v^{(i)}_{r+1}(U_i),\ldots, \rho^{(i)}_n X_n- v^{(i)}_n(U_i))_{\rho}$$

holds. Without loss of generality we may assume that $\Gamma_i$ represents
the coefficients of the polynomials $q_i$ and
$v_{r+1}^{(i)},\ldots,v_n^{(i)}$ with respect to $U_i$.

\end{thm}

\begin{pf} Under our hypotheses, namely that \xon are in Noether position with
respect to the variety $V_i$, the variables $\xo {n-i}$ being free, we have
the following integral ring extension of reduced rings:

$$R_i:=\Q[X_1,\ldots,X_{n-i}] \longrightarrow
B_i:=\Q[X_1,\ldots,X_n]/(f_1,\ldots,f_i).$$

Let $P_i =(p_1, \ldots,p_{n-i}) \in \Z^{n-i}$ be a lifting point of the
morphism $\pi_i$ and let $V_{P_i} = \pi_i^{-1}(P_i)$ be its lifting fiber.
We have $deg\ V_{P_i}=D_i=rank_{R_i} B_i \le deg\ V_i =\delta_i$. Let $U_i=
\lambda_{n-i+1}X_{n-i+1} + \cdots + \lambda_nX_n$$\in\Z[X_{n-i+1}\klk X_n]$
generate a primitive element of the ring extension $\Q\longrightarrow
\Q[V_{P_i}]$ (i.e. $U_i$ separates the points of $V_{P_i}$). The minimal
polynomial of the image of $U_i$ in $B_i$ has degree at least the
cardinality of the set $U_i(V_{P_i})$. Since the linear form $U_i$ separates
the points of $V_{P_i}$, this cardinality is $deg\ V_{P_i}= rank_{R_i} B_i$.
We conclude that $U_i$ generates also a primitive element of the ring
extension $R_i\longrightarrow B_i$. \spar

In the sequel we denote the primitive element generated by the linear form
$U_i$ in the ring extensions $\Q\longrightarrow\Q[V_{P_i}]$ and 
$R_i\longrightarrow B_i$ by the same letter $u_i$.\spar

By hypothesis  the following data are given explicitly (i.e. by the bit
representation of their coefficients):

\begin{itemize}

\item the primitive minimal equation, say $q\in\Z[T]$ of the primitive
element $u_i$ of $\Q[V_{P_i}]$,

\item the parametrization of $V_{P_i}$ by the zeroes of $q$, given by the
equations 

$$X_1-p_1=0, \ldots,X_{n-i}-p_{n-i}=0, q(T) = 0,$$ 
$$\rho_{n-i+1}X_{n-i+1}-v_{n-i+1}(T)=0, \ldots, \rho_n X_n- v_n(T)=0,$$

\end{itemize}

where $v_{n-i+1}\klk v_n$ are polynomials of $\Z[T]$ of degree strictly less
than $D_i=deg\ q$ and $\rho_{n-i+1}\klk \rho_n$ are non-zero integers. Check
again Definition \ref{geosolve} to see that the polynomials
$\rho_{n-i+1}X_{n-i+1}-v_{n-i+1}(T),\ldots, \rho_n X_n- v_n(T)$ are
assumed to be primitive.\spar

Let $a \in \Z$ be the leading coefficient of $q$.\spar

We consider now $f_1,\ldots,f_i$ as polynomials in the variables
$X_{n-i+1}\klk X_n$, i.e. as elements of the polynomial ring
$R_i[X_{n-i+1},\ldots,X_n]$. Let $f:=\ifo i $ and let

$$D(f):= \left( {{\partial f_k} \over {\partial X_j}} \right)_{ 1 \le k \le i
\atop n-i+1 \le j \le n}$$

be the corresponding
Jacobian matrix of $f$ (with respect to the variables $X_{n-i+1}\klk X_n$).
Recall that the Newton operator with respect to these 
variables is defined as

$$N_f =\pmatrix{X_{n-i+1}\cr \vdots\cr X_n \cr} - D(f)^{-1}
\pmatrix{f_1\cr \vdots\cr f_i \cr}.$$

Let $\kappa$ be a natural number. Using Lemma \ref{lemanewton} we deduce the
existence of numerators $g_{n-i+1}, \ldots,g_n$ and a non-zero denominator
$h$ in the polynomial ring $R_i[X_{n-i+1},\ldots,X_n]$ such that the
following $\kappa$-fold iterated Newton operator has the form:

$$N_f^{\kappa}= \pmatrix{{{g_{n-i+1}} \over {h}}\cr \vdots\cr{{g_n} \over
{h}} \cr}.$$

Now, let $M_{X_{n-i+1}}, \ldots,M_{X_n}$ be the matrices describing the
multiplication tensor of the $\Q$-algebra $\Q[V_{P_i}]$ (recall that by
assumption a geometric solution of the polynomial equation system defining
$V_{P_i}$ is given and that therefore the matrices $M_{X_{n-i+1}}\klk
M_{X_n}$ are known). 
\spar

Let $M$ denote the companion matrix of the polynomial $a^{-1}q(T) \in
\Q[T]$.  Let $n-i+1\leq j \leq n$. The following identity is immediate:

$$\rho_j M_{X_j}= v_j(M).$$

Moreover the matrices $\rho_j a^{D_i -1}M_{X_j}$ have integer entries. Let
$\kappa:=1+log_2\delta_i$ and note that $\kappa\geq 1+ log_2 D_i$ holds. Our
straight--line program $\Gamma_i$ will execute $\kappa$ Newton steps in a
subroutine which we are going to explain now: Let us consider the following
column vector of matrices $N_{n-i+1}\klk N_n$ with entries in
$\Q(X_1,\ldots,X_{n-i})$:

$$ \pmatrix{N_{n-i+1}\cr \vdots\cr N_n \cr} := N_f^{\kappa} ({\underline M})= 
 \pmatrix{{{g_{n-i+1}(X_1,\ldots,X_{n-i}, {\underline M})} \over
{h}(X_1,\ldots,X_{n-i}, {\underline M})}\cr \vdots\cr{{g_n(X_1,\ldots,X_{n-i},
{\underline M})} \over {h(X_1,\ldots,X_{n-i}, {\underline M})}} \cr},$$

where ${\underline M}= (M_{X_{n-i+1}}, \ldots,M_{X_n})$ and $g_{n-i+1}\klk
g_n$ are the numerator and $h$ the denominator polynomial of Lemma
\ref{lemanewton}. Finally, let us consider the matrix

$${\cal M}:= U_i(N_{n-i+1},\ldots,N_n)= \lambda_{n-i+1}N_{n-i+1}+ \cdots
+\lambda_nN_n.$$

This matrix is a matrix whose entries are rational functions of $\Q(\xo
{n-i})$.  From the fact that $P_i=(p_1\klk p_{n-i})$ is a lifting point and
from the proof of Lemma \ref{lemanewton} one deduces easily that in fact the
entries of ${\cal M}$ belong to the local ring

$$R_{P_i} := (R_i)_{(X_1-p_1 \klk X_{n-i}-p_{n-i})} =
\Q[X_1,\ldots,X_{n-i}]_{(X_1-p_1,\ldots,X_{n-i}-p_{n-i})}.$$

Let $T$ be a new variable. With these notations and assumptions we have the
following result:

\begin{lem}

Let $\chi \in \Q(X_1,\ldots,X_{n-i})[T]$ be the characteristic polynomial of
${\cal M}$ and $m_{u_i} \in R_i[T]$ the minimal integral equation of the
primitive element $u_i$ of $B_i$ over $R_i$.  Let $\chi(T) = T^{D_i}+
\sum\limits_{k=0}^{D_i -1} a_kT^k$ and $m_{u_i} = T^{D_i}+
\sum\limits_{k=0}^{D_i -1}b_k T^k$ with $a_k, b_k\in\Q(\xo {n-i})$ for
$0\leq k\leq D_i-1$. Then all these coefficients satisfy the condition
$ord_{P_i}(a_k-b_k) \ge \delta_i+1,$ where $ ord_{P_i}$ denotes the usual
order function (additive valuation) of the local ring $R_{P_i}$. 
\end{lem}

\begin{pf} Let $V_{P_i}= \{\xi_1,\ldots,\xi_{D_i} \}$ and fix $1\le l \le D_i$.
Let $\zeta_l\!:=U_i(\xi_l)$. Because of the hypothesis made on the lifting
fiber $V_{P_i}$ and from Hensel's Lemma (which represents a symbolic version
of the Implicit Function Theorem) we deduce that there exist formal power
series $R_{n-i+1}^{(l)},\ldots,R_n^{(l)} \in
\C[[X_1-p_1,\ldots,X_{n-i}-p_{n-i}]]$ with $R_{n-i+1}^{(l)}(P_i)
=\xi_{n-i+1}^{(l)},\ldots,R_n^{(l)}(P_i)=\xi_n^{(l)}$ such that for
$R^{(l)}\!:= (X_1-p_1,\ldots,X_{n-i}-p_{n-i}, R_{n-i+1}^{(l)},
\ldots,R_{n}^{(l)})$ the identities

\begin{equation} \label{fres0}
f_1(R^{(l)})=0,\ldots, f_i(R^{(l)})=0 
\end{equation} 

hold in $\C[[X_1-p_1,\ldots,X_{n-i}-p_{n-i}]]$. Let $u^{(l)}\!:=
U_i(R^{(l)}) = \lambda_{n-i+1}R_{n-i+1}^{(l)}\plp \lambda_nR_n^{(l)}$. As
shown in \cite{gihemomorpar} the minimal polynomial $m_{u_i}$ of the
primitive element $u_i$ verifies in
$\C[[X_1-p_1,\ldots,X_{n-i}-p_{n-i}]][T]$ the identity

\begin{equation}\label{minimo} m_{u_i}=\prod\limits_{1\le l \le  D_i}
(T-u^{(l)}).  \end{equation} 

Let us now consider the construction above which produces the matrix ${\cal
M} \in (R_{P_i})^{ D_i \times D_i}$ starting with the matrix $M \in \Q^{D_i
\times D_i}$ (recall that $M$ is the companion matrix of the polynomial
$a^{-1}q(T)\in\Q[T]$). The same construction transforms the $D_i$ distinct
eigenvalues of the diagonalizable matrix $M$, namely the values
$\zeta_l\!=U_i(\xi_l)$, $1 \le l \le D_i$, into eigenvalues of ${\cal M}$.
(Observe that by construction ${\cal M}$ is a rational function of the
matrix $M$ and therefore the same rational function applied to any
eigenvalue of $M$ produces an eigenvalue of ${\cal M}$). As shown in
\cite{gihemomorpar}, in this way we obtain $D_i$ distinct rational functions
${\tilde u}^{(l)} \in \C(X_1,\ldots,X_{n-i})$, which are eigenvalues of
${\cal M}$. Moreover these rational functions are all defined in the point
$P_i$ (this means that ${\tilde u}^{(l)} \in\C[X_1\klk
X_{n-i}]_{(X_1-p_1\klk X_{n-i}-p_{n-i})}$ holds). For $1\leq l\leq D_i$ the
rational functions ${\tilde u}^{(l)}$ can therefore be interpreted as
elements of the power series ring $\C[[X_1-p_1,\ldots,X_{n-i}-p_{n-i}]]$ and
they satisfy in this ring the condition

\begin{equation}\label{eqNr}
u^{(l)}-{\tilde u}^{(l)} \in
(X_1-p_1,\ldots,X_{n-i}-p_{n-i})^{\delta_i+1}.
\end{equation} 

(For a proof of these congruence relations see \cite{gihemomorpar}).
\spar

Let us now consider the characteristic polynomial $\chi$ of the matrix
${\cal M}$. Since the coefficients of ${\cal M}$ belong to $R_{P_i} =
\Q[X_1\klk X_{n-i}]_{(X_1-p_1,\ldots,X_{n-i}-p_{n-i})}$, the coefficients of
$\chi$ do too. Therefore $\chi$ can be interpreted as a polynomial in the
variable $T$ with coefficients in the power series ring $\Q[[X_1-p_1\klk
X_{n-i}-p_{n-i}]]$. From the fact that the rational functions ${\tilde
u}^{(1)}\klk {\tilde u}^{(D_i)}$ represent $D_i$ distinct eigenvalues of
${\cal M}$ we deduce that

\begin{equation} \label{chi} \chi(T) = \prod_{1\leq l\leq D_i} (T- {\tilde
u}^{(l)}) \end{equation}

holds. Let $\sigma_k$ denote the $k$--th elementary symmetric function in
$D_i$ arguments. The identity (\ref{chi}) implies that for $0\leq k \leq
D_i-1$ we can write the coefficient $a_k$ of $\chi(T)$ as
$a_k=(-1)^{D_i - k} \sigma_k({\tilde u}^{(1)},\ldots,{\tilde u}^{(D_i)})$.
>From the identity (\ref{minimo}) we deduce that the $k$-th coefficient $b_k$ of
the polynomial $m_{u_i}$ satisfies $b_k=(-1)^{D_i - k} \sigma_k(
u^{(1)},\ldots,u^{(D_i)})$. From the congruence relations (\ref{eqNr}) and
the identities (\ref{minimo}) and (\ref{chi}) we conclude now that for $0\leq
k\leq D_i-1$ the congruence relations

\begin{equation}\label{eqNoNo}
a_k-b_k \in  (X_1-p_1,\ldots,X_{n-i}-p_{n-i})^{\delta_i+1}
\end{equation}

hold in $\Q[[X_1-p_1,\ldots,X_{n-i}-p_{n-i}]]$. This proves the lemma.
\qed
\end{pf}

We continue now the proof of Theorem \ref{newton}. In order to see how to evaluate the
coefficients of $\chi$, let us consider $\theta:= det(h({\underline M})) \in
R_i$ and the matrix

$${\cal M}_1:= \theta {\cal M}.$$

This matrix ${\cal M}_1$ has entries in $R_i$. Now we have the identities

$$det(T \cdot Id_{D_i} - {\cal M})= det  (T \cdot Id_{D_i} - \theta^{-1} {\cal
M}_1) = \theta^{-D_i} det ((\theta T) Id_{D_i} - {\cal M}_1).$$

Let $\phi(T)=T^{D_i}+\phi_{D_i-1}T^{D_i-1}+ \cdots+ \phi_0 \in R_i[T]$
be the characteristic polynomial of ${\cal M}_1$. From the identities above
we deduce that for $0\leq k \leq D_i -1$ the $k$-th coefficient of $\chi$ 
can be written as

$$a_k= {{\theta^k \Phi_k} \over {\theta ^{D_i}}}.$$

Executing $\kappa = 1+log_2 \delta_i$ steps in the Newton iteration at the beginning
of this proof to produce the entries of the matrix $M$, $\theta$, ${\cal M}$
and finally ${\cal M}_1$ (in this order) and applying Lemma \ref{berk} we
produce a straight--line program $\Gamma_i'$ in $\Q(X_1\klk X_{n-i})$ of
non-scalar size $O(log_2\delta_i d^3i^6 L)$ and non-scalar depth $O((log_2 i
+ \ell)log_2\delta_i)$ using as parameters those given by the statement of
Theorem \ref{newton} such that $\Gamma_i'$ evaluates the family of
polynomials $1, \theta,\ldots, \theta^{D_i}$, $\phi_0,\ldots,\phi_{D_i-1}
\in R_i$. \spar

Now, applying Strassen's Vermeidung von Divisionen technique (Proposition
23) we obtain a division-free straight--line program $\Gamma_i''$ in $\Qxo
{n-i}$ of size $(id\delta_iL)^{O(1)}$, non-scalar depth $O((log_2 i +
\ell)log_2\delta_i)$ using as parameters the coordinates $p_1\klk p_{n-i}$
of $P_i$, the coefficients of the linear change of coordinates for the
Noether normalization for $V_i$, the rational numbers appearing as
coefficients in the geometric solution of the lifting fiber $V_{P_i}$ and
the parameters of $\Gamma$ such that $\Gamma_i''$ evaluates for each $0\leq
k\leq D_i -1$ the expansion of $a_k$ in $\Q[[X_1-p_1 \klk X_{n-i}-p_{n-i}]]$
up to terms of degree of order $\delta_i +1$. Taking into account the
congruence relations (\ref{eqNoNo}) we see that the division--free
straight--line program $\Gamma_i''$ evaluates polynomials $g_0,\ldots,
g_{D_i-1} \in \Q[X_1,\ldots,X_{n-i}]$ such that

$$ b_k -g_k \in (X_1-p_1,\ldots,X_{n-i}-p_{n-i})^{\delta_i+1}$$

holds in $\Q[[X_1-p_1\klk X_{n-i}-p_{n-i}]]$ for any $0\leq k\leq D_i-1$.
>From \cite{gihemomorpar} we deduce that the degrees of the coefficients
$b_k$ of the minimal polynomial $m_{u_i} \in R_i[T]$ do not exceed
$\delta_i$. Putting all this together, we conclude that

$$b_k = g_k$$

holds for any $0\leq k\leq D_i-1$. This means that the division-free
straight--line program $\Gamma_i''$ computes the coefficients of the
polynomial $q_i:= m_{u_i} \in \Q[\xo r, U_i] = R_i[U_i]$. \spar

In order to compute the parametrizations $v_{n-i+1}^{(i)}\klk v_n^{(i)}
\in\Q[\xo {n-i},T]$ we have to use the same kind of techniques (namely
truncated Newton--Hensel iteration) combined with the arguments developed in
\cite{krick-pardo1} and applied in \cite{gihemomorpar} and \cite{gihemopar}.
Let us be more exact: \spar

Let $n-i+1\leq k\leq n$ and let $m_{X_k}\in R_i[X_k]$ be the minimal
polynomial of the $R_i$-linear endomorphism given by the image of $X_k$ in
the $R_i$-algebra $B_i = \Q[V_i] = \Qxon / \ifo i$. The polynomial $m_{X_k}$
is monic and hence squarefree.  We compute the coefficients of the
polynomial $m_{X_k}\in R_i[X_k]$ in the same way as before the coefficients
of $q_i=m_{u_i}$. So, we have two monic squarefree polynomials $m_{X_k} \in
R_i[X_k]$ and $q_i \in R_i[U_i]$. Taking into account that $U_i$ separates
the associated primes of the ideal $(m_{X_k}, q_i)$ in $R_i[X_k,U_i]$, we
can apply directly \cite[Lemma 26]{krick-pardo1} in order to obtain the
parametrization associated to the variable $X_k$.  Doing the same for each
of the variables involved, namely the variables $X_{n-i+1},\ldots,X_n$, and
putting together the corresponding straight--line programs, we obtain a
procedure and a straight--line program $\Gamma_i$ of the desired complexity
which computes the output of Theorem 28. \qed
\end{pf}

\subsection{The recursion} \label{liftingpoint} 

\begin{prop}\label{punto} There exists a division-free arithmetic
network of non-scalar size of order $(i d \delta_i L)^{O(1)}$ and non-scalar depth
$O((log_2i + \ell)log_2 \delta_i)$ using parameters of logarithmic height
bounded by $max \{h, \eta_i, O((log_2i + \ell)log_2 \delta_i)\}$ which
takes as input

\begin{itemize}

\item a Noether normalization for the variety $V_i$, 

\item a lifting point $P_i$ and
 
\item a geometric solution of the lifting fiber $V_{P_i}$
\end{itemize}

and produces as output 

\begin{itemize} 

\item a linear change of variables $(X_1,\ldots,X_n) \longrightarrow
(Y_1^{(i+1)},\ldots,Y_n^{(i+1)})$ such that the new variables
$Y_1^{(i+1)},\ldots,Y_n^{(i+1)}$ are in Noether position with respect to
$V_{i+1}$,

\item a lifting point $P_{i+1}$ for $V_{i+1}$ and

\item a geometric solution for the lifting fiber $V_{P_{i+1}}$.

\end{itemize}

\end{prop} 

\begin{pf} 
The construction of the arithmetic network proceeds in three stages. In the
first stage we apply the algorithm underlying Theorem 28. In the second
stage we intersect algorithmically the variety $V_i$ with the hypersurfcae
$V(f_{i+1})$ in order to produce first a Noether normalization of the
variables with respect to the variety $V_{i+1} = V_i \cap V(f_{i+1})$. Then
we produce a linear form $U_{i+1}$ representing a primitive element
$u_{i+1}$ of the integral ring extension $R_{i+1} \longrightarrow B_{i+1}$
and a straight--line program representing polynomials analogous to the one
in the conclusion of Theorem 28. For this purpose we use the algorithm
underlying the proof of Proposition 14 in \cite{gihemomorpar}.
The only point we have to take care of is that we are working now over the
ground field $\Q$ and that we have to take into account the heights of the
parameters of $\Q$ we introduce in this straight--line program.
\spar

We remark that the
straight--line program in question has
non-scalar size $(id\delta_i L)^{O(1)}$ and non-scalar depth
$O(log_2(d\delta_i)+\ell)$ with parameters of logarithmic height bounded by
$O(log_2(d\delta_i) + \ell)$.
In the third and final stage we consider the polynomial

$$J(f_1\ldots,f_{i+1}):= det \left( {{\partial f_k}
\over {\partial X_j}} \right)_{ 1 \le k \le i+1 \atop n-i \le j \le n},$$

which is a nonzero divisor modulo the ideal 
$I_{i+1}=(f_1,\ldots,f_{i+1})$. Let us
observe that the polynomial $J(f_1\ldots,f_{i+1})$ can be evaluated by a
division-free straight--line program of length $O((i+1)^5+L)$ and depth
$O(log_2 (i+1) +\ell)$.
\spar
 
Let $\mu \in \K[X_1,\ldots,X_{n-i-1}]$ be the constant term of the
characteristic polynomial of the homothety given by $J(f_1\ldots,f_{i+1})$
modulo $I_{i+1}$. Since $J$ represents a nonzero divisor modulo $I_{i+1}$ we
conclude that $\mu$ does not vanish. Furthermore we observe that $\mu$ can be
evaluated by a division-free straight--line program in
$\Q[X_1,\ldots,X_{n-i-1}]$ of length $(i\delta_i L)^{O(1)}$ and depth
$O(log_2 d\delta_i + (log_2 i + \ell) log_2 \delta_i)= O((log_2i +
\ell)log_2 \delta_i)$, and so does the product $\rho \cdot \mu$, where
$\rho= \prod\limits_{n-i+1} ^n \rho_j^{(i)}$ is defined as in the statement
of Theorem 28.  \spar

Using a correct test sequence (see \cite{heschnorr}) we are able to find in
sequential time $(i\delta L)^{O(1)}$ and parallel time $O((log_2i +
\ell)log_2\delta_i)$ a rational point $P_{i+1} \in \Z^{n-i-1}$ of
logarithmic height bounded by $O((log_2i + \ell)log_2\delta_i)$ which
satisfies $(\rho \cdot \mu) (P_{i+1}) \not= 0$. Clearly, $P_{i+1}$ is a
lifting point for the variety $V_{i+1}$. In order to obtain a geometric
solution of $V_{P_{i+1}}$ with primitive element $u_{i+1}$ induced by the
linear form $U_{i+1}$ we have to specialize in the point $P_{i+1}$ the
polynomials obtained as output of the second stage (see \cite[Section
3]{gihemomorpar}). By this specialization we obtain the binary
representation of the coefficients of certain univariate polynomials in
$U_{i+1}$ which represent a geometric solution of the fiber $V_{P_{i+1}}$.
Nevertheless it might happen that the height of these coefficients is
excessive. In order to control the height of these polynomials we make them
primitive. This requires some integer greatest common divisor (gcd)
computations which do not modify the asymptotic time complexity of our
algorithm. \qed \end{pf}

\end{section} 
 

\typeout{Section 4} 

\begin{section}{Lifting Residues and Division modulo a
Complete Intersection Ideal}\label{division}

This section is dedicated to the proof of Theorem \ref{paso-division}.  The
outcome is a new trace formula for Gorenstein algebras given by complete
intersection ideals. This trace formula makes no reference anymore to a
given monomial basis of the algebra. Our trace formula represents an
expression which is ``easy-to-evaluate".


\begin{subsection}{Trace and Duality} 

Trace formulas appear in several recent papers treating problems in
algorithmic elimination theory. Some of these papers use a trace formula in
order to compute a quotient appearing as the result of a division of a given
polynomial modulo a given complete intersection ideal 
(see \cite{figismi,krick-pardo1}). Other papers use trace formulas in order to design
algorithms for geometric (or algebraic) solving of zero-dimensional
Gorenstein algebras given by complete intersection ideals 
(\cite{royetal,BeCaRoSz,Cardinal}). The paper \cite{saso} uses a trace formula to
obtain an upper bound for the degrees in the Nullstellensatz.  \spar

However, all these applications of trace formulas require the use of some
generating family of monomials of bounded degree which generate the given
Gorenstein algebra as a vector space over a suitable field. As a
consequence, such trace formulas provide just syntactical complexity or
degree bounds and in particular no intrinsic upper complexity bound (as e.g.
in Theorem \ref{paso-division}) can be obtained in this way. In this
subsection we introduce an alternative trace formula in order to get a
maximum benefit from the geometrically and algebraically well suited
features of Gorenstein algebras.  Let us start with a sketch of the trace
theory.  For proofs we refer to \cite{kunz}, Appendices E and F. \spar

Let $R$ be a ring of polynomials over a given ground field (for our
discussion the ground field may be assumed to be $\Q$).  Let $\K$ be the quotient
field of $R$ and let $R[X_1,\ldots,X_n]$ be the ring of $n$-variate polynomials
with coefficients in $R$. Let $f_1,\ldots,f_n$ be a smooth regular sequence
of polynomials in the ring $\Rxon$ of degree at most $d$ in the variables
$\xon$ generating a radical ideal denoted by $(f_1,\ldots,f_n)$. \spar

Consider now the $R$-algebra $B$ given as the quotient of $R[X_1,\ldots,X_n]$
by this ideal: 

$$B:=R[X_1,\ldots,X_n]/(f_1,\ldots,f_n).$$ 

We assume that the morphism $R\rightarrow B$ is an integral ring extension
representing a Noether normalization of the variety $V\ifon$ defined by the
polynomials $\fon$ in a suitable affine space. Thus, $B$ is a free
$R$-module of rank bounded by the degree of the variety $V(f_1,\ldots,f_n)$
(this estimation is very coarse but sufficient for our purpose).

Moreover, the $R$-algebra $B$ is Gorenstein and the following statements are
based on this fact. \spar

We consider $B^*:=Hom_R(B,R)$ as a $B$-module by the scalar product

$$ B\times B^* \longrightarrow B^*$$ 

which associates to any $(b,\tau)$ in $B\times B^*$ the $R$-linear map
$b\cdot \tau: B\longrightarrow R$ defined by $(b\cdot \tau) (x):= \tau(b x)$
for any element $x$ of $B$. \spar

Since the $R$-algebra $B$ is Gorenstein, its dual $B^{*}$ is a free
$B$-module of rank one.  Any element $\sigma$ of $B^{*}$ which generates
$B^{*}$ as $B$-module is called a {\sl trace} of $B$.  There are two
relevant elements of $B^*$ that we denote by $\Tr$ and $\sigma$. The first
one, \Tr, is called the {\em standard trace\/} of $B$ and it is defined in the
following way: given $b\in B$, let $\eta_b:B\longrightarrow B$ the
$R$-linear map defined by multiplying by $b$ any given element of $B$.  The
image $\Tr(b)$ under the map $\Tr$ is defined as the ordinary trace of the
endomorphism $\eta_b$ of $B$ (note that this definition makes sense since
$B$ is a free $R$-module). In order to introduce $\sigma$ (which will be a
trace of $B$ in the above sense), we need some additional notations. For any
element $g\in\Rxon$ we denote by $\bar{g}$ its image in $B$, i.e. the
residue class of $g$ modulo the ideal $\ifon$. Let $Y_1\klk Y_n$ be new
variables and let $Y:=(Y_1\klk Y_n)$. Let $1\leq j\leq n$ and let
$f_j^{Y}:=f_j(Y_1\klk Y_n)$ be the polynomial of $ \Ryon$ obtained by
substituting in $f_j$ the variables $\xon$ by $\yon$.  Let us consider the
polynomial

$$f_j^{Y}-f_j = \sum_{k=1}^n l_{jk}(Y_k-X_k) \in \R[\xon,\yon],$$

where the $l_{jk}$ are polynomials belonging to $R[\xon,\yon]$ having total
degree at most $(d-1)$ (observe that the $l_{jk}$ are not uniquely
determined by the sequence $\fon$). Let us now consider the determinant
$\Delta$ of the matrix $(l_{jk})_{1\leq j,k\leq n}$ which can be written
(non uniquely) as

$$ \Delta = \sum_m a_m(\xon) b_m(\yon) \in \R[\xon,\yon],$$

with the $a_m$ being elements of $\Rxon$ and the $b_m$ elements of $\Ryon$.
\hyphenation{deter-minant} (Observe that it will not be necessary to find
the polynomials $a_m$ and $b_m$ alegbraically, we need just their existence
for our argumentation.) The polynomial $\Delta$ is called a pseudo-jacobian
determinant of the regular sequence \ifon. Observe that the polynomials
$a_m$ and $b_m$ can (and will) be chosen to have degrees bounded by
{\hbox{$n(d-1)$}} in the variables $\xon$ and $\yon$ respectively. Let
$c_m\in\Rxon$ be the polynomial we obtain by substituting in $b_m$ the
variables $\yon$ by $\xon$. For $\bar{J}$ the class of the Jacobian
determinant $J\ifon$ in $B$ we have the identity

$$\bar{J} = \sum_m \bar{a}_m\cdot\bar{c}_m.$$

Moreover the image of the polynomial $\Delta$ in the residue class ring
\linebreak 
$R[\xon,\yon]$ modulo the ideal $(\fon,f_1^Y\klk f_n^Y)$ is independent of
the particular choice of the matrix $(l_{kj})_{1\leq k,j\leq n}$. This
justifies the name ``pseudo-jacobian" for the polynomial $\Delta$.  With
these notations there is a unique trace $\sigma \in B^*$ such that the
following identity holds in $B$:

$$\bar 1 = \sum_m \sigma (\bar{a}_m)\cdot\bar{c}_m.$$

The main property of the trace $\sigma$, known as ``trace formula" (``Tate's
trace formula" \cite[Appendix F]{kunz}, \cite{Iversen} being a special case of it) is
the following statement: for any $g\in\Rxon$ the identity

\begin{equation}\label{eqTrNo}
\bar g = \sum_m\sigma(\bar g \cdot
\bar{a}_m)\cdot \bar{c}_m
\end{equation}

holds true in $B$. Let us observe that the polynomial $\sum_m\sigma(\bar g
\cdot \bar{a}_m)\cdot c_m \in \Rxon$ underlying the identity (\ref{eqTrNo})
is of degree at most $n(d-1)$ in the variables $\xon$.

The main use of this trace formula consists in solving the following
problem:

\begin{prob}[Lifting of a residual class] \label{problem6} Given a
polynomial \linebreak
$g\in\Rxon$ of arbitrary degree in $\xon$, find a polynomial
$g_1\in\Rxon$ of degree at most $n(d-1)$ in the variables $\xon$ such that
$\bar{g}_1 = \bar g$ holds in $B$. \end{prob}

As we have seen before, the trace formula (\ref{eqTrNo}) solves Problem
\ref{problem6} since it allows us to choose for $g_1$ the polynomial

\begin{equation}\label{eqTrNoNo} 
g_1 := \sum_m\sigma(\bar g \cdot \bar{a}_m)\cdot c_m.
\end{equation}

However, defining the polynomial $g_1$ by the formula (\ref{eqTrNoNo})
inhibits us from taking advantage of any special ``semantical" features of
the $R$-algebra $B$: one ``a priori" needs all monomials of degree at most
$n(d-1)$ for the description of the polynomials $c_m$ (and $a_m$).
Therefore, we replace the trace formula (\ref{eqTrNo}) by the following
alternative one:

\begin{prop}[Trace Formula] With the same notations as before, let us
consider the free \Rxon-module $\Bxon$ given by extending scalars in $B$
(this means we consider the tensor product $\Bxon := B\otimes_R \Rxon$) and
let us also consider the polynomial $\Delta_1\in\Rxon$ given by:

$$\Delta_1:=\sum_m\bar{a}_m\cdot c_m \in \Bxon.$$

Then for any $g\in\Rxon$ the following identity holds true in \linebreak
$\Rxon$:

$$\sum_m\sigma(\bar g \cdot\bar{a}_m)\cdot c_m =
\Trt (\bar{J}^{-1} \bar g \cdot \Delta_1).$$

(Here $\Trt:= \Tr \otimes Id_{\Rxon} : \Bxon\longrightarrow \Rxon$ is the
standard trace obtained from the standard trace $\Tr: B\longrightarrow R$ by
extending scalars).

\end{prop}

\begin{pf} Let $\Tr: B\longrightarrow R$ be the standard trace of the free
$R$-module $B$. Let us recall from \cite[Appendix F]{kunz} that for any
$g\in\Rxon$ the identity

$$\Tr(\bar{J}^{-1}\bar g) = \sigma(\bar g)$$

holds. From the $\Ryon$-linearity of the map $\Trt:\Byon \longrightarrow
\Ryon$ we deduce that any $g\in\Rxon$ satisfies the identities

$$\Trt(\bar{J}^{-1}\bar g \Delta_1(\yon)) = 
\sum_m \Trt(\bar{J}^{-1}\bar g \bar{a}_m\cdot b_m) = 
\sum_m \Trt(\bar{J}^{-1}\bar g \bar{a}_m) \cdot b_m.$$ 

In other words we have in $\Ryon$

$$\Trt(\bar{J}^{-1}\bar g \Delta_1(\yon)) = 
\sum_m \sigma(\bar g \cdot \bar{a}_m) \cdot b_m$$ 

for any $g\in\Rxon$. 
Replacing in this identity the variables $\yon$ by $\xon$ we obtain the
desired formula
 
$$\Trt(\bar{J}^{-1}\bar g \Delta_1) = 
\sum_m \sigma(\bar g \cdot a_m) \cdot c_m\ .\mbox{\hskip 3cm \qed}$$
\end{pf}

One sees easily that for any $h\in R[\xon,\yon]$,
$\Trt(\bar h)$ is simply the standard trace of the image $\bar h$ of $h$ in
the $\Ryon$-module $\Byon$. This observation together with Proposition 31
represents our basic tool for the evaluation of formula (\ref{eqTrNoNo}) and
hence for the solution of Problem \ref{problem6}. This is the content of the
following observations.
\spar

Let us consider the $\K$-algebra $B' = \K\otimes_R B$ obtained by localizing
$B$ in the non-zero elements of $R$. Fix a basis of the finite dimensional
$\K$-vectorspace $B'$. Let $M_{X_1} \klk M_{X_n}$ be the matrices of the
homotheties $\eta_{X_i}: B' \longrightarrow B'$ with respect to the given
basis of $B'$ and let $\Tr$ denote the function which associates to a given
matrix its usual trace. With these conventions let $g_1$ be defined as

\begin{eqnarray}\label{eqlast41}
g_1 & := & \Tr( J(\fon) {(M_{X_1} \klk M_{X_n})}^{-1} 
           \cdot g(M_{X_1} \klk M_{X_n}) \nonumber \\
    &    & \cdot \Delta(M_{X_1} \klk M_{X_n}, \xon) ) .
\end{eqnarray}

One easily verifies that $g_1$ belongs to $\Rxon$ and that $\bar g_1 = \bar
g$ holds in $B$. 

\end{subsection} 


\begin{subsection}{A Division Step} 

The lifting process presented in the last subsection is now applied to compute
the quotient of two polynomials modulo a reduced complete intersection
ideal. More precisely, let us consider $f\in\Rxon$ a polynomial which is not
a zero-divisor in $B$ and another polynomial $g\in\Rxon$ such that the
residual class $\bar f$ divides the residual class $\bar g$ in $B$.  The
following proposition shows how we can compute a lifting quotient $q\in\Rxon$
for the division of $\bar g$ by $\bar f$ in $B$.

\begin{prop}[Division Step]\label{Prop32} Let notations and assumptions be
the same as in the previous subsection and let $D$ be the rank of $B$ as
free $R$-module. Let be given the following items as input:

\begin{itemize}

\item a straight--line program $\Gamma'$ of size $L$ and depth $\ell$
representing the polynomials $f, g$ and $\fon$.

\item the matrices $M_{X_1}\klk M_{X_n}$ describing the multiplication
tensor of $B$ with respect to the given basis of $B'=\K\otimes_R B$.

\end{itemize}

Suppose that $\bar f$ is a non-zero divisor of $\bar B$ and that $\bar f$
divides $\bar g$ in $B$. Then there exists a division-free straight--line
program $\Gamma$ in $\Kxon$ of size $L(n d D)^{O(1)}$ and non-scalar depth
$O( \ell + log_2 D + log_2n)$ which computes from the entries of the
matrices $M_{X_1}\klk M_{X_n}$ and the parameters of $\Gamma'$ a non-zero
element ${\theta}$ of $R$, and a polynomial $q$ of $\Rxon$ such that
$\theta$ divides $q$ in $\Rxon$ and such that $\bar q \bar f = \bar \theta
\bar g$ holds in $B$.

\end{prop}

\begin{pf} In order to prove this result, let us observe that any basis of $B$ as
free $R$-module induces a basis of $\Byon$ as free \Ryon-module. Moreover if
$M_{X_i}$ is the matrix of the multiplication by $\bar{X}_i$ in $B$ with
respect to a given basis, $M_{X_i}$ represents as well the multiplication by
$\bar{X}_i$ in $\Byon$ with respect to the same basis. Next, since the
polynomials $f$ and $J\ifon$ are not zero-divisors modulo \ifon, the
following matrices are non-singular:

$$ F_1 : = f(M_{X_1} \klk M_{X_n}),$$

$$J_1:=det ( {{\partial f_i} \over {\partial X_j}}(M_{X_1} \klk
M_{X_n}))_{1\leq i,j\leq n}.$$

Finally, let us denote by $G_1$ and $\Delta_1$ the following two matrices:

$$G_1:= g(M_{X_1} \klk M_{X_n}),$$

and

$$\Delta_1:=\Delta(M_{X_1} \klk M_{X_n},\yon),$$

where $\Delta$ is the pseudojacobian determinant of $\fon$. Let us remark
that the matrices $F_1,J_1$ and $G_1$ have entries in $\K$ while $\Delta_1$
has entries in\linebreak 
$\K[\yon]$. From formula (\ref{eqlast41}) of the previous
subsection we deduce that $q_1 := \Tr(J_1^{-1}\cdot F_1^{-1}\cdot G_1\cdot
\Delta_1(\xon))$ is a polynomial of $\Rxon$ which satisfies in $B$ the
identity $\bar{q}_1\bar f =\bar g$ in $B$ (by $\Tr$ we denote here the
ususal trace of matrices).

Finally, let us transpose the adjoint matrices of $F_1$ and $J_1$:

$$F_2 := \;\;^t\!Adj(F_1),\mbox{ and } J_2:= \;\;^t\!Adj(J_1).$$

The  quotient $q\in\Rxon$ and the non-zero constant ${\theta} \in R$
we are looking for are given as follows:

\begin{itemize}
\item $q:=\Tr(F_2 \cdot J_2 \cdot G_1\cdot \Delta_1)$
\item ${\theta} := det(F_2)\cdot det(J_2)$
\end{itemize}

Clearly, the polynomial $q_1$ is the quotient of $q \over {\theta}$, while
$\bar q \bar f = \bar{\theta} \bar g$ holds in $B$.  Moreover $q$ can be
computed by a {\em division-free\/} straight--line program $\Gamma$ in $\Kxon$
from the entries of $M_{X_1}\klk M_{X_n}$ and the parameters of $\Gamma'$
(note that for the computation of $q_1$ we need divisions). The complexity
bounds in the statement of Proposition \ref{Prop32} follow by reconstruction of the
straight--line programs that evaluate $f, g, J\ifon, \Delta$ and the
determinants involved (cf. Subsection \ref{elementary} above). \spar

{\sl Proof of Theorem \ref{paso-division} }
\spar

In order to prove Theorem \ref{paso-division}  we just follow the algorithm
underlying the proof of Proposition \ref{Prop32}. We have to add just some comments
concerning the matrices $M_{X_1}\klk M_{X_n}$. Let us consider a linear form
$u = \lambda_1 X_1 \plp \lambda_nX_n$ (with $\lambda_i\in\Z$) inducing a
primitive element of the zero-dimensional $\Q$-algebra $\Qxon / \ifon$. In
this case $R$ will be the field $\Q$. Let $q_u\in\Z[T]$ be the minimal
polynomial of $u$ and $\rho_1X_1-v_1(T)\klk \rho_nX_n-v_n(T)$ the
parametrizations of the variety with respect to this primitive element. Let
$\alpha$ be the leading coefficient of $q_u$ and $\rho = \prod_{i=1}^n
\rho_i$ a discriminant. Then the companion matrix of $q_u$ has the form

$$\alpha^{-1} M,$$

where $M$ is a matrix with integer entries.  Now, the matrices 
describing the multiplication tensor of $\Qxon / \ifon$ can be written as

$$ M_{X_i} = \rho_i^{-1}\cdot v_i(\alpha^{-1} M)$$

for $1\leq i\leq n$. Taking $g=1$ and $f=f_{n+1}$ in Proposition
\ref{Prop32} we obtain ${\theta} \in \Q, {\theta}\neq 0$ and $q\in\Qxon$
such that

$$ {\theta} \cdot 1 - q \cdot f_{n+1} \in \ifon$$

holds. Finally, multiplying by appropriate powers of $\alpha$ and $\rho$
we obtain a non-zero integer $a$ and a polynomial $g_{n+1}$ of the form

$$ a:=\alpha^N\rho^M\cdot{\theta} \in \Z \ \ \ \ \ \ \ g_{n+1} :=
\alpha^N\rho^M q\in \Zxon,$$

such that $a-g_{n+1}f_{n+1}\in \ifon$ holds. The bounds of the Theorem
\ref{paso-division}
are then obtained from the bounds of Proposition \ref{Prop32}. The bounds for the
height of the parameters are obtained simply by choosing an appropriate
primitive element $u$ such that $q_u$ and the parametrizations have height
equal to the minimal height of the diophantine variety $V:=V\ifon$. \qed
\end{pf}

\end{subsection} 

\end{section} 

\begin{ack}
L. M. Pardo wants to thank the \'Ecole Polytechnique X for the invitation and 
hospitality during his stay in the Fall of 1995, when this paper was
conceived.
\end{ack}

\typeout{References}

\end{document}